\documentclass{amsart}
\usepackage{amssymb, amsbsy, amsthm, amsmath, amstext, amsopn}

\newtheorem{thm}{Theorem}[section]
\newtheorem{lemma}[thm]{Lemma}
\newtheorem{cor}[thm]{Corollary}
\newtheorem{prop}[thm]{Proposition}
\newtheorem*{extlemma}{Extension Lemma}
\newtheorem*{uniqlemma}{Uniqueness Lemma}

\theoremstyle{definition}

\newtheorem{defn}[thm]{Definition}

\newtheorem*{notation}{Notation}

\theoremstyle{remark}

\newtheorem{remark}[thm]{Remark}
                  
\newtheorem{claim}[thm]{Claim}
                
\newtheorem{const}[thm]{Construction}

\newtheorem*{pfsketch}{Sketch of Proof}

\usepackage[all]{xy}
\usepackage{amsfonts}
\usepackage{amscd}

\begin{document}

\newcommand{\hs}{\mbox{\hspace{.4em}}}
\newcommand{\ds}{\displaystyle}
\newcommand{\bd}{\begin{displaymath}}
\newcommand{\ed}{\end{displaymath}}
\newcommand{\bcd}{\begin{CD}}
\newcommand{\ecd}{\end{CD}}
\newcommand{\theo}{\mathcal{O}}
\newcommand{\proj}{\operatorname{Proj}}
\newcommand{\bproj}{\underline{\operatorname{Proj}}}
\newcommand{\spec}{\operatorname{Spec}}
\newcommand{\bspec}{\underline{\operatorname{Spec}}}
\newcommand{\pline}{{\mathbf P} ^1}
\newcommand{\boldp}{{\mathbf P}}
\newcommand{\boldq}{{\mathbf Q}}
\newcommand{\cs}{{\mathbf C} ^*}
\newcommand{\boldc}{{\mathbf C}}
\newcommand{\M}{{\mathcal M}}
\newcommand{\bM}{{\mathbf M}}
\newcommand{\ldb}{[[}
\newcommand{\rdb}{]]}
\newcommand{\Sym}{\operatorname{Sym}^{\bullet}}
\newcommand{\Symp}{\operatorname{Sym}}
\newcommand{\lrar}{\leadsto}
\newcommand{\cK}{{\mathcal K}}
\newcommand{\cB}{{\mathcal B}}
\newcommand{\cE}{{\mathcal E}}
\newcommand{\cA}{{\mathcal A}}
\newcommand{\cF}{{\mathcal F}}
\newcommand{\cI}{{\mathcal I}}
\newcommand{\cL}{{\mathcal L}}
\newcommand{\cQ}{{\mathcal Q}}
\newcommand{\cU}{{\mathcal U}}
\newcommand{\Aut}{\operatorname{Aut}}
\newcommand{\End}{\operatorname{End}}
\newcommand{\Hom}{\operatorname{Hom}}
\newcommand{\uHom}{\underline{\operatorname{Hom}}}
\newcommand{\Ext}{\operatorname{Ext}}
\newcommand{\bExt}{\operatorname{\bf{Ext}}}
\newcommand{\inv}{^{-1}}
\newcommand{\airtilde}{\widetilde{\hspace{.5em}}}
\newcommand{\airhat}{\widehat{\hspace{.5em}}}
\newcommand{\Hilb}{\operatorname{Hilb}}
\newcommand{\supp}{\operatorname{supp}}
\newcommand{\id}{\operatorname{id}}
\newcommand{\nt}{^{\circ}}

\newcommand{\psstack}{\M_S(L_1\oplus L_2)}
\newcommand{\TF}{\operatorname{TF}}
\newcommand{\Bun}{\operatorname{Bun}}

\title[Framed sheaves on ruled surfaces]{Moduli Spaces of Framed Sheaves on Certain Ruled Surfaces over Elliptic Curves}
\author{Thomas A. Nevins}
%% \date{\today}
\address{Department of Mathematics\\University of Michigan\\525 East University Avenue, Ann Arbor, MI 48109-1109 USA}
\email{nevins@math.lsa.umich.edu}

\begin{abstract}
Fix a ruled surface $S$ obtained as the projective completion of a line bundle $L$ on a complex elliptic curve $C$; we study the moduli problem of parametrizing certain pairs consisting of a sheaf $\cE$ on $S$ and a map of $\cE$ to a fixed reference sheaf on $S$.  We prove that the full moduli stack for this problem is representable by a scheme in some cases.  Moreover, the moduli stack admits an action by the group $\cs$, and we determine its fixed-point set, which leads in some special cases to explicit formulas for the rational homology of the moduli space.
\end{abstract}

\maketitle

\tableofcontents

\numberwithin{equation}{section}

\section{Introduction}

In this work, we compute the rational homology of moduli spaces (or stacks) of rank two framed 
torsion-free sheaves on certain complex ruled surfaces over elliptic curves.

Moduli spaces of framed sheaves on surfaces have previously been an object of study in a variety of
contexts (\cite{MR95i:14015}, \cite{MR96a:14017}, \cite{MR95k:58026}, \cite{MR1718147}); however, we believe that there are now good reasons for concentrating 
significant attention on a very special class of these moduli spaces, namely the spaces $\M_S(E)$ 
parametrizing
sheaves on the projectivization $S$ of the total space of a line bundle $L$ on a smooth projective 
curve $C$ that are equipped with isomorphisms of their restrictions to the divisor at infinity 
$D\subset S$ with a fixed reference sheaf $E$ on $D$.

One reason for focusing attention on this case
is that it provides a tool for computations relevant to the geometric Langlands program for 
surfaces.  Given a curve $C$ in a surface $S'$, one wants to compute the algebra of Hecke operators 
arising from modification of sheaves on $S'$ along the curve $C$.  If one allows only point 
modifications along $C$, then one knows the algebra completely, thanks
to work of Nakajima (\cite{MR97j:17027}, \cite{MR98h:14006}),
Grojnowski (\cite{MR97f:14041})
and Baranovsky (\cite{math.AG/9811092}).  If, however, one allows modification along the entire curve $C$, the computation of the algebra
appears to be extremely difficult, although progress in a special case by Nakajima (\cite{MR95i:53051},
\cite{MR99b:17033})  and Ginzburg--Kapranov--Vasserot (\cite{MR96f:11086}, \cite{KV}) suggests that the picture should be a very rich one.

Informed by the study of Donagi--Ein--Lazarsfeld in pure dimension 1 (\cite{MR98f:14006}), one might
try to obtain information concerning Hecke operators for $C\subset S'$ by degenerating the inclusion 
$C\subset S'$ to the inclusion of $C$ in the normal cone of $C$ in $S'$; the algebra of Hecke operators here then ought to be
a degeneration of the full algebra for $C\subset S'$.  This ``normal cone setting'' is the one 
developed in this paper: taking $L$ to be the normal bundle of a smooth curve $C$ in a smooth 
surface $S'$, the moduli space $\M_S(E)$ provides the appropriate setting for computing the 
degenerate algebra of Hecke operators.  We make a start at computing the rational homology of our 
moduli spaces (focusing on the special case in which $L$ is a degree zero line bundle over an elliptic curve) as a necessary preliminary to understanding this degenerate algebra; note that the
full algebra seems not to be known even when $C$ is an elliptic curve (but see \cite{MR96f:11086} for a very interesting
conjecture), and we expect our technique to give new information even in this simple case.  In another paper (\cite{framedlocalize}), we develop tools for studying these moduli spaces for much more general line bundles $L$ over curves of higher genus.

A variant of our moduli spaces $\M_S(E)$ arises also in connection with singular Higgs pairs on $C$ in the work
of Jardim (\cite{math.DG/9912028}, \cite{math.DG/9910120}, \cite{math.AG/9909146}, \cite{math.DG/9909069}),
and $\M_S(E)$ appears (in the special case $L=K_C$) to figure in the study of
$\mathcal{D}$-modules on curves, the adelic or Beilinson--Drinfeld Grassmannian, and the moduli of 
bundles on noncommutative surfaces (see \cite{MR1670188}, \cite{MR99f:58107}, and \cite{hep-th/0002193} for some related work that suggests this relationship).

In the present work, we concentrate on two rather concrete questions.  The first of these is the representability of the
moduli stack $\M_S(E)$ by a scheme or algebraic space.  This study is motivated in part by a desire
to study further some more classically geometric features of these stacks (for example, the natural
$L^2$-metrics and existence of integrable system structures) and in part by a desire to deploy 
a powerful tool---localization---in the study of the homology of $\M_S(E)$.  The stack $\M_S(E)$ has
a $\cs$ action coming from the action on $S$, and we wish to use this action to answer our second
question: how to produce homology
bases for the stack.  However, to apply localization, one
generally needs a smooth separated algebraic manifold (or smooth Hausdorff K\"{a}hler manifold) for 
which limits $\ds\lim_{\lambda\rightarrow 0}\lambda\cdot m$ exist for all points $m$ of the manifold (here $\lambda\in\cs$).  This requirement on the limits of orbits for the $\cs$ action makes it problematic to apply the usual stability or semistability conditions to obtain good moduli spaces that retain all the features necessary for the application of localization techniques, and this leads to some
interplay between our answers to the two questions.

In Section \ref{stacks} we define the moduli stack $\M_S(E)$ in question and describe
carefully the $\cs$ action, and in Section \ref{fixedpoints} we use this to determine the fixed points of the
action on $\M_S(E)$.  In Section \ref{chapterextensionlemma} we develop the tools that show the existence of the limits
we need and the properness of the components of the fixed-point set of $\cs$ in $\M_S(E)$ and, also, with some additional work, allow us to prove separatedness of $\M_S(E)$  in some cases.
Then in Section \ref{chrepresentability} we prove our principal representability and smoothness results, focusing on the case of an elliptic curve $C$.  In the related paper \cite{framedlocalize}, we will add to these representability results
a resolution theorem that allows one to disregard issues of separatedness and stackiness of the
moduli problem in all genera; this resolution theorem opens the door for application of the singular localization
theorems of \cite{MR89c:14072} and \cite{MR85d:32063} in a very general setting.
Finally, in Section \ref{chhomologybasis} we use the work of earlier sections to compute a basis in rational homology for a particular moduli space $\psstack$ when the curve $C$ is elliptic and $S\rightarrow C$ satisfies an additional (open) condition.

\vspace{0.07in}

The author is deeply indebted to his dissertation advisor, Kevin Corlette, 
without whose help, guidance, and encouragement this work could not have
been completed.  He also wishes to thank Madhav Nori for many discussions, in particular for his suggestion that the author's earlier, much more complicated statement and proof of Corollary \ref{noricor} might be explained by something like Proposition \ref{norissuggestion}; Brendan Hassett and R. Narasimhan for helpful advice; and Vladimir Baranovsky, Stanley Chang, Victor Ginzburg, Michael Mandell, Stephanie Nevins, Tony Pantev, Amritanshu Prasad, and Ian Robertson for helpful conversations.  The author's graduate work at the University of Chicago, of which this paper is a result, was supported in part by an NDSEG fellowship from the Office of Naval Research.

\section{The Moduli Stack: Definition and Group Action}\label{stacks}

In this section, we define the moduli stacks in which we are interested and describe some relevant properties of the group action.

\subsection{Notation and Definition of the Moduli Stack}

Fix a smooth complete irreducible complex curve $C$, a line bundle $L$ on $C$, and a rank two vector bundle $E$ on $C$; later we will specialize to the case in which $C$ has genus one.  Let
\bd
S={\mathbf P}(L\oplus \theo),
\ed
and let $D$ denote the divisor at infinity in $S$.  We will use the notation $\sigma$ for the curve $C$ embedded in $S$ as the zero section of the bundle $L$. Since $D\cong C$ canonically via the projection map $\pi: S\rightarrow C$, we may think of $E$ as a vector bundle on $D$ or on $\sigma$.  If $R$ is a $\boldc$-scheme, let $E_R$ denote the pullback of $E$ to $D_R = D\times R$.

\begin{defn} 
Let $\M_S(E)$ denote the moduli stack of $E$-framed rank two torsion-free sheaves on $S$; this stack has as its objects pairs $(\cE_R,\phi_R)$ consisting of an $R$-flat family of rank two torsion-free sheaves $\cE_R$ on $S_R$ together with an isomorphism
\bd
\phi_R: \cE_R\big|_{D_R} \xrightarrow{\sim} E_R.
\ed
A morphism $(\cE_R,\phi_R) \rightarrow (\cE'_T,\phi'_T)$ consists of a morphism $f: R\rightarrow T$ together with an isomorphism 
\bd
\begin{CD}
\cE_R@>\psi>\sim> (1_S\times f)^*\cE'_T
\end{CD}
\ed
for which 
\bd
\xymatrix{\cE_R\big|_{D_R}\ar[d]_{\psi} \ar[rr]^{\phi_R} & & E_R\ar[d]^{\operatorname{id}}\\
(1_S\times f)^*\cE'_T\big|_{D_R} \ar[rr]^-{(1_D\times f)^*\phi_T} & & E_R}
\ed
commutes.
\end{defn}

\begin{prop}(\cite{framedlocalize})
$\M_S(E)$ is an algebraic stack that is locally of finite type over $\spec\boldc$.
\end{prop}
\begin{pfsketch}
Let $\TF_S$ denote the moduli stack parametrizing torsion-free sheaves on $S$, and let $\TF_S(D)$ denote its open substack that parametrizes those sheaves that are locally free along $D$ (see \cite{LMB}).  There is a restriction morphism
\bd
\TF_S(D) \rightarrow \Bun(D),
\ed
where $\Bun(D)$ denotes the moduli stack parametrizing vector bundles on $D$.  The vector bundle $E$ on $D$ determines a morphism $\spec\boldc \rightarrow \Bun(D)$, and $\M_S(E)$ is the fiber product of the diagram
\bd
\xymatrix{ & \TF_S(D) \ar[d] \\
\spec\boldc \ar[r] & \Bun(D).}
\ed
The proposition follows because $\spec\boldc\rightarrow\Bun(D)$ is an $\Aut(E)$-torsor over a locally closed substack of $\Bun(D)$.\,\,\qedsymbol
\end{pfsketch}

\subsection{The Group Action on the Surface $S$}

Let $L$ denote a line bundle on the curve $C$.  Let $\theo$ denote the trivial bundle on $C$, which we think of as having the global generator $s$.  Let
\bd
S= \bproj \Sym (L^* \oplus \theo^*).
\ed

Define an action of $\cs$ on $\Sym (L^*\oplus \theo^*)$ by
\begin{equation}\label{normalaction}\begin{split}
\lambda\cdot l^* & = l^*,\\
\lambda\cdot s^* & = \lambda\inv s^*.
\end{split}\end{equation}
Notice that this induces the action on the underlying space of
\mbox{$\bproj\Sym(L^*\oplus\theo^*)$} that extends the action of $\cs$ on $\bspec\Sym L^*$, although
 its action on the ring $\Sym (L^*\oplus\theo^*)$ does not extend that on $\Sym L^*$.

Recall that an action of $\cs$ on $\pline$ is given, in more geometric
terms, by a map
\bd
\cs\times\pline \xrightarrow{m} \pline
\ed
and hence is determined by the map of rings
\bd
\theo_{\pline} \longrightarrow m_*\theo_{\cs\times\pline} =
\theo_{\pline}[\lambda,\lambda\inv].
\ed
Globally on $S$, then, we want a map
\bd
\Sym (L^*\oplus \theo^*) \xrightarrow{\overline{m}} \Sym
(L^*\oplus\theo^*)[\lambda,\lambda\inv];
\ed
to obtain the action in Equation \eqref{normalaction} as the appropriate
specialization at a fixed $\lambda\in\cs$, we must take
\begin{equation}\begin{split}
L^* \mapsto L^*,\\
\theo^* \mapsto \lambda\inv\theo^*.
\end{split}\end{equation}

\begin{remark}
In our applications to $\M_S(E)$, we will want to pull sheaves back
along the {\em inverse} action, i.e., the action by which
$\lambda\in\cs$ acts on $S$ via the usual action of $\lambda\inv$ on
$S$.  It is convenient in our discussion to compute everything for the
standard action on $S$ and then later to simply reverse the weights of
the $\cs$ action when we are concerned with the inverse action.
\end{remark}

Choose a local generator $l^*$ of $L^*$.  We will translate between $\Sym(L^*\oplus\theo^*)$ and $\boldc [l^*,s^*]$ and thus describe some of the usual structures on $S$.

The zero section of $L\subset S$ is just the set along which the section $l^*$ vanishes, hence
\bd
\theo_\sigma = \Sym (L^*\oplus\theo^*)/\langle L^*\rangle,
\ed
where $\langle L^*\rangle$ denotes the ideal generated by $L^*$.  Now $\langle L^*\rangle$ is exactly the image of 
\bd
L^*\otimes \Sym (L^*\oplus \theo^*)
\ed
in $\Sym(L^*\oplus\theo^*)$ under multiplication, so
\bd
\theo (-\sigma)=\bigl[ L^*(-1)\otimes \Sym (L^*\oplus\theo^*)\bigr]\airtilde.
\ed
Here the notation $L^*(-1)$ indicates that the grading on $L^*$ is shifted so that it now lies in graded degree 1; this makes the inclusion
\bd
L^*(-1)\otimes\Sym (L^*\oplus\theo^*) \subset\Sym (L^*\oplus\theo^*)
\ed
compatible with the gradings and consequently gives a map of the associated coherent sheaves (indicated by the tilde) on $\bproj\Sym (L^*\oplus\theo^*)$.

This description determines $\theo (k\sigma)$ for $k\in\mathbf{Z}$: the sheaf $\theo (\sigma)$ is the dual of $\theo (-\sigma)$ or
\bd
\theo(\sigma) = \bigl[ L(1)\otimes\Sym(L^*\oplus\theo^*)\bigr]\airtilde,
\ed
and $\theo (k\sigma)$ is therefore
\bd
\theo(k\sigma) = \bigl[ L^k(k)\otimes \Sym (L^*\oplus\theo^*)\bigr]\airtilde.
\ed
In terms of the coordinate ring $\boldc [l^*,s^*]$ we get
\bd
\theo (k\sigma) = (l^*)^{-k}\boldc [l^*,s^*]\airtilde,
\ed
where here $(l^*)^{-k}$ lies, as expected, in graded degree $-k$. 

Recall that a structure of $\cs$-equivariant $\theo_S$-module on
$\widetilde{M}$ is given by an isomorphism
\begin{equation}
m^*M\xrightarrow{I} M[\lambda,\lambda\inv]
\end{equation}
of $\Sym(L^*\oplus\theo^*)[\lambda,\lambda\inv]$-modules with certain
compatibility conditions (see \cite{MR97j:14001}).  The action of  an element 
$\lambda\in\cs$
on sections of $\widetilde{M}$ is then given by the composition
\begin{equation}
M\xrightarrow{m^*} m^*M \xrightarrow{I}
M[\lambda,\lambda\inv]\xrightarrow{\operatorname{eval}(\lambda)} M.
\end{equation}
We want to use the equivariant structures on submodules of
$\theo_S$ induced from the canonical equivariant structure on $\theo_S$ that comes from the group action on $S$.

 Over $S\setminus D$, the sheaf $\theo (k\sigma)$ consists of
 elements of $(l^*)^{-k}\boldc [l^*,s^*]_{s^*}$ of degree zero; this
 module has generator $\ds\left( \frac{s^*}{l^*}\right) ^k$, on which
 $\cs$ acts with weight $-k$ owing to the $(s^*)^k$ term in the
 numerator: under the map
\begin{equation}
\theo_S \longrightarrow \theo_S[\lambda,\lambda\inv],
\end{equation}
the generator $\ds\left( \frac{s^*}{l^*}\right) ^k$ of $\theo(k\sigma)$ (for $k\leq 0$) maps to
$\ds\lambda^{-k}\left(\frac{s^*}{l^*}\right)^k$, so this element lies
in the weight $-k$ subspace.  

\begin{remark}
Note, again, that this implies that the generator lies in the weight
$k$ subspace for the {\em inverse} action.  We will summarize our
results for the inverse action in Proposition \ref{actionsummary}.
\end{remark}

Over $S\setminus\sigma$ we find that $\theo (k\sigma)$ consists of elements of $(l^*)^{-k}\boldc [l^*,s^*]_{l^*}$ of degree zero, and this module has as its generator the unit of the ring $\boldc [s^*/l^*]$, on which $\cs$ acts with trivial weight.  Thus, $\cs$ acts with weight zero in the fibers of $\theo (k\sigma)$ along $D$ but with weight $-k$ in the fibers of $\theo (k\sigma)$ along $\sigma$.  

It is convenient here to note also that, since the ring $\theo_S$ on $S\setminus \sigma$ comes from $\boldc [s^*/l^*]$, the group $\cs$ acts on functions defined near $D$ as follows: if $f$ is a function homogeneous in the variable $s^*/l^*$, then $\cs$ acts on $f$ with weight equal to minus the order of vanishing of $f$ along $D$.     Similarly, near $\sigma$ the group $\cs$ acts with weight $k$ on the functions homogeneous in the variable $l^*/s^*$ vanishing to order exactly $k$ along $\sigma$.

Recall that the space of global sections of a module $\theo (n) = \widetilde{\Gamma(\theo (n))}$ over $\pline$ is just the space of degree zero elements of the graded module $\Gamma_* (\theo (n))$ (see \cite{MR57:3116} for this notation).  The construction of our line bundles allows us to conclude that the modules by means of which we defined these line bundles are already in the image of the functor $\Gamma_*$, and hence, if
\bd
\pi : S\rightarrow C
\ed
is the canonical projection, that
\bd
\pi_* \theo (k\sigma) = \left(L^k\otimes\operatorname{Sym}^k (L^*\oplus \theo^*)\right).
\ed
Of course if $k<0$ this is the zero sheaf; however, if $k\geq 0$ one obtains
\begin{align}
\pi_*\theo (k\sigma) & = \sum_{j=0}^k L^k\otimes (L^*)^j\otimes (\theo^*)^{k-j}\\
& = \sum_{j=0}^k L^{k-j}\otimes (\theo^*)^{k-j}.
\end{align}
If we identify this direct sum with 
\begin{equation}
\pi_*\theo (k\sigma) = \sum_{j=0}^k L^{k-j},
\end{equation}
then our earlier description implies that $\cs$ acts on the factor $L^{k-j}$ with weight $j-k$.

We want to include one further, similar calculation.  Performing the same analysis as above for
\bd
\theo (-D) = \langle \theo^*\rangle \subset\Sym(L^*\oplus\theo^*),
\ed
one gets 
\bd
\theo (-kD) = (\theo^*)^k(-k)\otimes \Sym(L^*\oplus\theo^*).
\ed
In terms of the ring $\boldc [l^*,s^*]$, one then obtains
\bd
(s^*)^k(-k)\boldc [l^*,s^*],
\ed
and one may describe this sheaf near $D$ as the sheaf associated to the module of degree zero elements of 
\bd
(s^*)^k\boldc [l^*,s^*]_{l^*}.
\ed
This module has generator $\left(s^*/l^*\right)^k$ and so $\cs$ acts on its generator with weight $-k$.  In particular, we obtain the following conclusion.

\begin{lemma}
Let $I_D$ denote the ideal sheaf of $D$.  Then $\cs$ acts on
\bd
\left( N^*_{C/S}\right)^k = \frac{I^k_D}{I^{k+1}_D} \cong L^k
\ed
with weight $-k$.
\end{lemma}

\begin{remark}
The description in this section of weights of elements in various sheaves easily allows computation of weights in, for example, $\cs$-invariant finite-colength ideal sheaves or their quotients as well, simply because these weights are induced from those in the structure sheaf of $S$.
\end{remark}

We summarize our results for the inverse action.

\begin{prop}\label{actionsummary}
For the inverse action
\bd
m\inv (\lambda, s)= \lambda\inv\cdot s
\ed
of $\cs$ on $s\in S$ and the induced equivariant structures on sheaves coming from $\theo_S$,
\begin{enumerate}
\item the homogeneous generator of $\theo (k\sigma)$ along $\sigma$ lies in weight $k$, and the homogeneous generator along $D$ lies in weight zero;
\item $\ds\pi_*\theo(k\sigma) = \sum_{j=0}^k L^{k-j}$\newline and $\cs$ acts on $L^{k-j}$ with weight $k-j$;
\item $\cs$ acts on $I_D^k/I_D^{k+1} \cong L^k$ with weight $k$; and
\item $\cs$ acts on the homogeneous functions that vanish to order $k$ along $\sigma$ with weight $-k$.
\end{enumerate}
\end{prop}

\subsection{The Group Action on $\M_S(E)$}

\begin{const}
The stack $\M_S(E)$ admits a $\cs$ action as follows: $\lambda\in\cs$ acts by pulling back along the map $m_{\lambda^{-1}}:S\rightarrow S$ given by multiplication by $\lambda^{-1}$ on $S$ fiberwise, which extends the action on the total space of the line bundle $L$ by scaling in the fibers.
\end{const}

Suppose that
\bd
0\rightarrow M_1 \rightarrow M_2 \rightarrow M_3 \rightarrow 0
\ed
is a short exact sequence of $\theo_S$-modules and that $M_1$ and $M_3$ are equipped with $\cs$-equivariant structures.  Suppose $M_2$ is defined by a \v{C}ech $1$-cocycle $\theta \in Z^1\bigl(S, \uHom (M_3,M_1)\bigr)$.  Then $m_{\lambda\inv}^*M_2$ occurs as a 1-extension
\bd
0\rightarrow m_{\lambda\inv}^*M_1 \rightarrow m_{\lambda\inv}^*M_2 \rightarrow m_{\lambda\inv}^* M_3 \rightarrow 0,
\ed
and, if
\begin{align}
I_1: m_{\lambda\inv}^*M_1 & \longrightarrow M_1,\\
I_3: m_{\lambda\inv}^*M_3 & \longrightarrow M_3
\end{align}
are the isomorphisms coming from the $\cs$-equivariant structure, then $m_{\lambda\inv}^*\theta$ yields the $1$-cocycle 
\bd
I_1 \circ m_{\lambda\inv}^*\theta \circ I_3\inv \in Z^1 \bigl( S, \uHom (M_3,M_1)\bigr)
\ed
that determines $m_{\lambda\inv}^*M_2$ as an extension of $M_3$ by $M_1$.  Hence $m_{\lambda\inv}^*M_2$ is determined by the class 
\bd
\bigl[ \lambda\inv \cdot \theta\bigr] \in H^1\bigl(S,\uHom (M_3,M_1)\bigr).
\ed

In order to determine the invariant $1$-extensions, then, we need only determine the invariants in such cohomology (or, more generally, Ext) groups.  The machinery of the next section makes it possible to use this data to determine fixed points of $\M_S(E)$.

\section{The Fixed-Point Set}\label{fixedpoints}

We begin by describing a technical device that allows concrete study of sheaves on $S$ by means of extensions.

\subsection{The Canonical Short Exact Sequence}

Let $S=\boldp (L\oplus \theo)$, $L$ a line bundle on the curve $C$.
We will describe the canonical short exact sequence associated to a
rank two torsion-free sheaf on $S$.

\subsubsection{The Short Exact Sequence for Sheaves of Type {\em U}}\label{exactseqtypeU}
We first consider the case of a rank two torsion-free sheaf on $S$ of
{\em unequal generic splitting type}, or {\em type U}; by this we mean a sheaf $\cE$ whose restriction to a generic fiber $f$ of $S\rightarrow C$ is of the form
\bd
\cE\big| _f \cong \theo (d_1)\oplus \theo (d_2)
\ed
with $d_1\neq d_2$.  Suppose $\cE$ is such a sheaf that is in addition locally free, and suppose $d_1>d_2$.  The {\em canonical short exact sequence} associated to $\cE$ is the sequence
\begin{equation}
0\rightarrow \bigl[\pi^*\left(\pi_*\cE(-d_1\sigma)\right)\bigr](d_1\sigma) \rightarrow \cE \rightarrow Q\rightarrow 0.
\end{equation}
This exact sequence has been used widely to study sheaves on rational and ruled surfaces; we will use the papers of Brosius (\cite{MR85i:14009},
\cite{MR85i:14010}) as our canonical references for this sequence.

\begin{prop}(Brosius)
If $\cE$ is locally free of type $U$, then $Q$ is a rank one torsion-free sheaf
on $S$; hence the exact sequence is a 1-extension of the form
\bd
0\rightarrow L_1 \rightarrow \cE \rightarrow L_2\otimes I_2 \rightarrow 0
\ed
where $L_1$ and $L_2$ are line bundles on $S$ and $I_2$ is an ideal of finite colength.  If $\cE_T$ is a $T$-flat family of rank two vector bundles on $S$ where $T$ is some integral Noetherian base scheme and the generic splitting type 
\bd
\cE_t\big|_f \cong \theo(d_1)\oplus \theo(d_2)
\ed
is constant as $t\in T$ varies, then the canonical subsheaf $(L_1)_T$ for the family specializes properly:
\bd
((L_1)_T)_t = (L_1)_t.
\ed
\end{prop}

There is a similar short exact sequence of the form
\begin{equation}
0\rightarrow L_1\otimes I_1\rightarrow \cE\rightarrow L_2\otimes I_2\rightarrow 0
\end{equation}
when $\cE$ fails to be locally free; here $L_1$ and $L_2$ are line bundles and $I_1$ and $I_2$ are finite colength ideals.  This short exact sequence has two useful descriptions.
\begin{enumerate}
\item\label{firstdescription} Given $\cE$ torsion-free of generic splitting type 
\bd
\cE\big| _f = \theo (d_1)\oplus \theo (d_2),
\ed
let 
\bd
\mathcal{L} := \left[ \pi^*\left(\pi_*\cE(-d_1\sigma)\right)\right](d_1\sigma).
\ed
Then $L_1\otimes I_1\subset \cE$ is the saturation of $\mathcal{L}$ in $\cE$, and by definition the quotient $\cE/(L_1\otimes I_1)$ is thus a torsion-free sheaf of rank one.

\item Apply the construction of Brosius to $\cE^{**}$; one gets a diagram
\begin{displaymath}
\xymatrix{0\ar[r] & L_1' \ar[r] & \cE^{**} \ar[r] & L_2'\otimes I_2' \ar[r] & 0 \\
0 \ar[r] & L_1'\cap \cE \ar[r]\ar[u] & \cE \ar[r]\ar[u] & L_2\otimes I_2 \ar[r]\ar[u] & 0}
\end{displaymath}
for which the vertical arrows are inclusions.  Then $L_1\otimes I_1 = L_1'\cap \cE$ is exactly as in description \ref{firstdescription}.
\end{enumerate}

\begin{remark}
Notice that if $\cE$ is $\cs$-equivariant, then by construction both
$\bigl(\pi^*\left[\pi_*\cE(-d_1\sigma)\right]\bigr)(d_1\sigma)$ and (as then follows from the first description) $L_1\otimes I_1$ are $\cs$-invariant subsheaves of $\cE$.
\end{remark}

One may show, using techniques similar to those employed by Brosius, that the canonical exact sequence specializes properly when $\cE$ varies in families of torsion-free sheaves, provided one fixes the proper invariants; however, we will not need this fact.

\subsubsection{The Short Exact Sequence for Sheaves of Type {\em E}}

Suppose now that $\cE$ is of {\em equal generic splitting type} or
{\em type E}: this means that the restriction of $\cE$ to the generic
fiber $f$ is of the form
\bd
\cE\big|_f \cong \theo (d)\oplus \theo (d).
\ed
We will assume that $\cE$ is locally free, since this will suffice for
our intended application later in the section.  In this case, one
obtains a canonical short exact sequence
\begin{equation}
0\rightarrow \bigl(\pi^*E'\bigr)\otimes\theo(d\sigma) \rightarrow \cE
\rightarrow \mathcal{I}_{Z\subset X}\otimes \theo (d\sigma) \rightarrow 0,
\end{equation}
where 
\bd
\bigl(\pi^*E'\bigr)\otimes\theo(d\sigma) =
\pi^*\bigl(\pi_*\cE(-d\sigma)\bigr) (d\sigma)
\ed
and $\mathcal{I}_{Z\subset X}$ is the ideal associated to the
inclusion of an l.c.i. (that is, local complete intersection) zero-cycle $Z\subset S$ in the scheme
\bd
X= \pi\inv \left(\pi(Z)\right) \subset S.
\ed
In fact, if $\pi^*E'$ is the vector bundle above, then $Z$ is exactly
the zero-cycle along which the map
\bd
\pi^*E' \rightarrow \cE(-d\sigma)
\ed
of rank two vector bundles fails to be a rank two linear map of the fibers.

\subsection{The Fixed Points for the Action}

We proceed to study the fixed points of $\cs$ in $\M_S(E)$.  We begin by proving that any simple fixed point admits a $\cs$-equivariant structure.

In what follows, let
\bd
m: \cs\times S \rightarrow S
\ed
denote the (inverse) multiplication map given by
\bd
m(\lambda, s) = \lambda ^{-1} \cdot s.
\ed
We also let $m_{\lambda}$ denote the map $m(\lambda, \cdot):S\rightarrow S$ given by multiplication by $\lambda\inv$.

\begin{lemma}\label{uniqueequivariance}
Let $R$ denote any $\boldc$-scheme.  Suppose $\cE_R$ is an $R$-flat family of $R$-simple torsion-free rank two coherent sheaves on $S$ that is equipped with a
framing
\bd
\cE_R\big|_{D_R} \xrightarrow{\phi_R} E_R.
\ed
Suppose furthermore that for each $\lambda\in\cs$, one has
\bd
m_\lambda ^*(\cE_R,\phi_R) \cong (\cE_R,\phi_R).
\ed
Then $\cE_R$ admits a $\cs$-equivariant structure for which $\phi_R$ is $\cs$-equivariant.
\end{lemma}
\begin{proof}
Let 
\bd
p:\cs\times S\rightarrow S
\ed
denote the projection on the second factor, let
\bd
q:S\times R\rightarrow R
\ed
denote the projection to $R$, and, by abuse of notation, let
\bd
q: D\times R \rightarrow R
\ed
also denote the projection to $R$.

By assumption, we have
\bd
m_\lambda^*\cE_R \cong \cE_R
\ed
for each $\lambda\in\cs$.  Since $\cE_R$ is $R$-simple, we find that
\bd
\Ext_{q|\lambda}^0 \left(m_\lambda ^*\cE_R,(p^*\cE_R)_\lambda\right)
\ed
is a line bundle on $R$ for each $\lambda\in\cs$; therefore, we get a line bundle 
\bd
\Ext_{(\id\times q)}^0 (m^*\cE_R,p^*\cE_R)
\ed
on $\cs\times R$.  There is a restriction map
\bd
\Ext_{(\id\times q)}^0 (m^*\cE_R,p^*\cE_R) \rightarrow \Ext_{(\id\times q)}^0(m^*\cE_R\big| _{D_R},p^*\cE_R\big| _{D_R});
\ed
since
\bd
m^*\cE_R\big|_{D_R} \cong E_{\cs\times R} \text{\hspace{.5em}and\hspace{.5em}} p^*\cE_R\big|_{D_R}\cong  E_{\cs\times R}
\ed
canonically via $\phi_R$, we find that
\begin{align*}
\Ext_{(\id\times q)}^0(m^*\cE_R\big|_{D_R},p^*\cE_R\big|_{D_R}) & \cong \Ext_{(\id\times q)}^0 (E_{\cs\times R},E_{\cs\times R})\\
& \cong \theo_{\cs\times R}
\end{align*}
canonically.  Because each $m_\lambda^*(\cE_R,\phi_R)$ is isomorphic to $(\cE_R,\phi_R)$, the induced map of line bundles
\bd
\Ext_{(\id\times q)}^0 (m^*\cE_R,p^*\cE_R) \rightarrow \theo_{\cs\times R}
\ed
is surjective and therefore is an isomorphism (see Theorem 2.4 of \cite{MR90i:13001}).  The inverse image of the identity section $1\in\theo_{\cs\times R}$ now gives a section of $\Ext_{(\id\times q)}^0 (m^*\cE_R,p^*\cE_R)$; moreover, the unicity of our construction guarantees commutativity of the necessary diagram (see Section 3.2 of \cite{MR97j:14001}) to make the pulled-back morphism
\bd
m^*\cE_R \rightarrow p^*\cE_R
\ed
an equivariant structure on $\cE_R$.  This completes the proof.
\end{proof}

\begin{remark}
In the concrete cases in which we will be interested in Sections \ref{chrepresentability} and \ref{chhomologybasis}, all $E$-framed sheaves that are invariant under each $\lambda\in\cs$ will admit $\cs$-equivariant structures.  In fact, the proper notion of a $\cs$-fixed point in our stack is a $\cs$-equivariant framed sheaf in general; however, because we will not need this level of generality in this work, we will omit a detailed explanation of the reasons that this is the appropriate ``fixed-point substack.''
\end{remark}

\subsubsection{Fixed Points: The Unequal Case}

Now suppose that $(\cE,\phi)$ defines a $\cs$-invariant point of $\M_S(E)(\spec\boldc)$ and that $\cE$ is of type {\em U}. Replacing $\cE$ by its double dual, we assume for the moment that $\cE$ is a locally free and $\cs$-equivariant rank two coherent sheaf on $S$.  The short exact sequence
\bd
0\rightarrow I_1\otimes L_1 \rightarrow \cE \rightarrow I_2\otimes L_2 \rightarrow 0
\ed
is then, by construction, also $\cs$-equivariant.  Since $\cE$ is locally free, we must have $I_1=\theo_S$.

\begin{lemma}\label{notalongD}
The support of $\theo/I_2$ has empty intersection with $D$.
\end{lemma}

\begin{proof}
Because $\cE\big| _D \cong E$ by a $\cs$-equivariant map, we may find, locally along $D$, some 
$\cs$-invariant sections $s_1$ and $s_2$ of $\cE$ that generate $\cE$.  These sections induce a $\cs$-invariant isomorphism $\cE \cong \theo ^2$.  Thus, locally along $D$ we get a $\cs$-invariant map $L_1 \rightarrow \theo^2$.  Choosing a local $\cs$-invariant generator of $L_1$ gives an expression of this map as
\bd
\psi: \theo \rightarrow \theo^2,
\ed
say with $\psi(s) = (f_1\cdot s, f_2\cdot s)$.  Because $\psi$ is $\cs$-equivariant, we find that both $f_1$ and $f_2$ are $\cs$-invariant.  But now if the image of $(f_1,f_2)$ vanishes in a fiber $\cE_p$ of $\cE$ at a point $p\in D$, then $(f_1,f_2)$ vanishes along $\pi^{-1}(\pi(p))$, hence $\cE/(I_1\otimes L_1)$ has torsion supported along $\pi^{-1}(\pi(p))$.  However, this quotient is torsion-free by construction, and we may conclude that $(f_1,f_2)$ fails to vanish at any point $p\in D$.  Therefore $\cE/(I_1\otimes L_1)$ is locally free at every point $p\in D$, completing the proof.
\end{proof}

Next, assume, without loss of generality, that 
\bd
\cE\big|_f \cong \theo(l\sigma)\oplus \theo((-1-l)\sigma)
\ed
along a general fiber $f$, where $l\geq 0$.  Write 
\bd
L_1 = \theo_S(l\sigma)\otimes \pi^*B_1
\ed
and
\bd
L_2 = \theo_S((-1-l)\sigma)\otimes \pi^*B_2.
\ed

\begin{lemma}\label{i2istrivial} 
$I_2 = \theo_S$.
\end{lemma}

\begin{proof}
We have an exact sequence
\begin{multline*}
\Ext^1\left( (\theo/I_2)\otimes L_2,L_1\right)^{\cs} \rightarrow \Ext^1(L_2,L_1)^{\cs}\rightarrow \\
\Ext^1(I_2\otimes L_2,L_1)^{\cs} \rightarrow \Ext^2\left( (\theo/I_2)\otimes L_2,L_1\right)^{\cs}.
\end{multline*}
Now $\Ext^1\left( (\theo/I_2)\otimes L_2,L_1\right) = 0$, since it is (after applying Serre duality) isomorphic to $H^1$ of a sheaf supported on a finite collection of points.  Consequently, it is enough to show that $\Ext^2\left( (\theo/I_2)\otimes L_2,L_1\right)^{\cs} = 0$: then any equivariant extension pulls back from an extension of $L_2$ by $L_1$, which can only be locally free if already $L_2 = I_2\otimes L_2$.  So, we must compute $\Ext^2\left( (\theo/I_2), L_2^{-1}\otimes L_1\right)^{\cs}$.  Since\bd
L_2^{-1}\otimes L_1 = \theo((1+2l)\sigma) \otimes \pi^*(B_2^{-1}\otimes B_1)
\ed
and $\theo/I_2$ has discrete support,
we may replace $L_2^{-1}\otimes L_1$ by $\theo((1+2l)\sigma)$ for this computation.  Now, locally along $\sigma$ the quotient $\theo/I_2$ has generators of the form $z^{k_1}(l^*/s^*)^{k_2}$ with $z$ a uniformizer in the $\sigma$ direction and $l^*/s^*$ a parameter in the fiber direction.  Here $\cs$ acts by
\bd
\lambda \cdot (z^{k_1}(l^*/s^*)^{k_2}) = \lambda^{-k_2}z^{k_1}(l^*/s^*)^{k_2};
\ed
hence, using short exact sequences in the first variable of $\Ext$, we may reduce to considering $\Ext^2\left(\boldc_{\chi}, \theo((1+2l)\sigma)\right)^{\cs}$, where $\cs$ acts on $\boldc_{\chi}$ with nonpositive weight $\chi$.  If we resolve $\boldc_{\chi}$ locally by
\bd
\theo(-f-\sigma)\underset{\boldc}{\otimes}\boldc_\chi  \rightarrow \left[\theo(-f)\oplus\theo(-\sigma)\right] \underset{\boldc}{\otimes} \boldc_\chi \rightarrow \theo\underset{\boldc}{\otimes} \boldc_\chi \rightarrow \boldc_\chi,
\ed
we find that $\Ext^2\left(\boldc_\chi, \theo ((1+2l)\sigma)\right)^{\cs}$ is an invariant quotient of 
\bd
H^0\left(\theo((2+2l)\sigma)\otimes \theo(f)\underset{\boldc}{\otimes}\boldc_{-\chi}\right)^{\cs}.
\ed
But any $\cs$-invariant section of
\bd
\theo((2+2l)\sigma)\underset{\boldc}{\otimes} \boldc_{-\chi}
\ed
vanishes to order at least $2+2l - \deg \chi$ along $\sigma$ and thus maps to zero in 
\newline $\Ext^2\left( \boldc_\chi, \theo((1+2l)\sigma)\right)$.\end{proof}

Consequently, a $\cs$-equivariant locally free sheaf is a $\cs$-invariant extension of the form
\bd
0\rightarrow L_1 \rightarrow \cE \rightarrow L_2\rightarrow 0.
\ed

\begin{lemma}\label{equalextgroups}
The map 
\bd
\Ext^1(L_2,L_1)^{\cs} \rightarrow \Ext^1(L_2,L_1\big|_D)^{\cs}
\ed
is an isomorphism.
\end{lemma}

\begin{proof}
\bd
\Ext^1(L_2,L_1) = H^1(S,L_2^{-1}\otimes L_1)
= H^1\left(S, \theo ((1+2l)\sigma)\otimes \pi^*(B_2\inv\otimes B_1)\right).
\ed
The Leray spectral sequence shows that this is just
\bd
H^1\left( D, \pi_*\theo ((1+2l)\sigma)\otimes B_2\inv\otimes B_1\right);
\ed
furthermore, the restriction map is evidently just the map
\bd
H^1\left(D, \pi_*\theo ((1+2l)\sigma)\otimes B_2\inv\otimes B_1\right) \rightarrow H^1\left( D, \theo_D\otimes B_2\inv\otimes B_1\right),
\ed
which is given fiberwise by the map $H^0\left( \pline, \theo ((1+2l)\cdot p)\right) \rightarrow \boldc$ (where $p$ is the zero point in $\pline$) arising from evaluation at infinity.  
Since the unique $\cs$-invariant section of $\theo_{\pline}((1+2l)\cdot p)$ under our normalization of the $\cs$ action is the one that vanishes to order $1+2l$ at $p$ (and hence is nonzero at infinity), the map
\bd
\left[\pi_*\theo_S ((1+2l)\sigma)\right]^{\cs} \rightarrow \theo_D
\ed
is an isomorphism, and the conclusion follows.
\end{proof}

Finally, then, returning to the case in which $\cE$ is not locally free, we see that $\cE^{**}/\cE$ is a finite-length $\cs$-invariant quotient of $\cE^{**}$.  We may summarize as follows.

\begin{prop}
Every $\cs$-equivariant $\cE$ is the kernel of a finite-length $\cs$-equivariant quotient of a locally
free $\cs$-equivariant coherent sheaf $\cF$.  The sheaf $\cF$ is uniquely determined by the sheaf $L_1$ in its canonical short exact sequence 
\begin{equation}\label{shortexact}
0\rightarrow L_1 \rightarrow \cF \rightarrow L_2 \rightarrow 0
\end{equation}
together with the inclusion $L_1\big|_D \subset E$.
\end{prop}

\begin{proof}
Given $\cE$, let $\cF = \cE^{**}$.  Then the extension in Equation \ref{shortexact} determines data of $L_1$ and the inclusion $L_1\big|_D \subset E$, and Lemma \ref{equalextgroups} shows that this data in turn determines $\cF$ as a framed bundle.  Now the inclusion $\cE\subseteq \cF$ is $\cs$-equivariant, hence the quotient $\cF/\cE$ is also $\cs$-equivariant, as desired.
\end{proof}

\subsubsection{Fixed Points: The Equal Case}

Suppose now that $\cE$ is locally free of type {\em E}.  Let 
\bd
0\rightarrow (\pi^*E')(d\sigma)\rightarrow \cE \rightarrow {\mathcal I}_{Z\subset Y}(d\sigma)\rightarrow 0
\ed
denote the canonical exact sequence.  The same argument as we used in Lemma \ref{notalongD} for the type U case shows that, if $\cE$ is $\cs$-invariant, then $\operatorname{supp} Z\subset\sigma$.  

Consider the boundary map
\bd
\Ext^1\bigl({\mathcal I}_{Z\subset Y}(d\sigma), (\pi^*E')(d\sigma)\bigr) \rightarrow \Ext^2\bigl(\theo_Z(d\sigma), (\pi^*E')(d\sigma)\bigr)
\ed
coming from the short exact sequence in the first variable
\bd
0\rightarrow {\mathcal I}_{Z\subset Y} \rightarrow \theo_Y \rightarrow \theo_Z\rightarrow 0.
\ed
As in the case of type U, any class in $\Ext^1\bigl({\mathcal I}_{Z\subset Y}(d\sigma), (\pi^*E')(d\sigma)\bigr)$ whose image in $\Ext^2\bigl(\theo_Z(d\sigma), (\pi^*E')(d\sigma)\bigr)$ is zero must define a 1-extension that is not locally free, unless $Z=\emptyset$.  But
\bd
\Ext^2\bigl(\theo_Z(d\sigma), (\pi^*E')(d\sigma)\bigr)^{\cs} \cong \bigl[\Ext^2(\theo_Z,\theo)^2\bigr]^{\cs}
=0
\ed
because all weights in $\Ext^2(\theo_Z,\theo)$ are greater than or equal to one.  So any $\cs$-invariant $1$-extension has singularities unless $Z=\emptyset$.

This analysis gives the following proposition.

\begin{prop}
Every $\cs$-invariant $E$-framed sheaf of type E is the kernel of a $\cs$-equivariant morphism
\bd
(\pi^*E)(d\sigma) \rightarrow {\mathcal Q},
\ed
where ${\mathcal Q}$ is a finite-length $\cs$-equivariant $\theo_S$-module supported along $\sigma$.
\end{prop}

\section{The Extension Lemma and its Relatives}\label{chapterextensionlemma}

We first introduce some notation for this section.  We fix a line bundle $B$ on a smooth complete curve $C$, and let $\cA$ denote the quasicoherent $\theo_C$-algebra 
\bd
\cA = \Sym_{\theo_C}B.
\ed
We introduce also several variations on this notation: we let $I\subset\cA$ denote the ideal $\ds\oplus_{n\geq 1} \operatorname{Sym}^n_{\theo_C}B$, $\widehat{\cA}$ denote the completion of $\cA$ with respect to $I$ (note that $\widehat{\cA}$ is no longer a quasicoherent sheaf!), and $\cA^\circ$ denote the localization of $\cA$ with respect to $I$.  If $T$ is a $\boldc$-scheme, we let $\widehat{(\cA_T)}$ denote the completion of the pullback of $\cA$ to $T$ and $\cA_T^\circ$ denote the pullback of $\cA^\circ$ to $T$.  Most often in this context we will let $R$ denote a DVR with uniformizing element $\pi$ and field of fractions $K$, and then we will use either $\spec R$ or $\spec K$ as our scheme $T$; in these cases, we will abbreviate to $\cA_R$, etc.

The following fundamental lemma, which is a simplified version of a result of Langton from \cite{MR51:510}, is essential in establishing separatedness of some moduli stacks and ``orbit completion'' properties.

\begin{extlemma}\label{extlemma} (\cite{framedlocalize}) Suppose $\cE_K$ is a torsion-free coherent $\cA$-module on $C_K$ of rank two, $\cE_R^\circ$ is an $R$-flat family of torsion-free coherent $\cA^\circ$-modules on $C_R$, and 
\bd
\psi: \cE_K \otimes_{\cA_K} \cA_K^\circ \longrightarrow \cE_R^\circ \otimes_{\cA_R^\circ} \cA_K^\circ
\ed
is an isomorphism of $\cA^\circ$-modules on $C_K$.  Then there is a unique $\cA_R$-submodule $\cE_R\subseteq \cE_K$ of $\cE_K$ for which
\begin{enumerate}
\item $\cE_R \big|_{C_K} = \cE_K$,
\item $\ds\cE_R\otimes_{\cA_R}\cA_R^\circ \big| _{\bspec \cA^\circ} \cong \cE_R^\circ \big| _{\bspec \cA^\circ}$ compatibly with $\psi$, and
\item $\cE_R$ is an $R$-flat family of torsion-free $\cA$-modules on $C_R$.
\end{enumerate}
\end{extlemma}

The proof of this result uses the equivalence of the categories of quasicoherent $\cA$-modules and quasicoherent sheaves on $\bspec\cA$ to translate between commutative algebra and geometry.

The Extension Lemma has the following important corollary.

\begin{cor}\label{corofextlemma}(\cite{framedlocalize}) Let $S={\mathbf P}(L\oplus \theo)$.  Let 
\bd
m: \cs \times S \longrightarrow S
\ed
denote the multiplication map given by
\bd
m(\lambda, s) = \lambda^{-1} \cdot s.
\ed
Fix a rank two torsion-free coherent sheaf $\cE$ on $S$, with restriction to the divisor at infinity $\cE\big|_{D} = E$ a rank two vector bundle $E$ on $D$.  Then $m^*\cE$ extends to a rank two flat family of torsion-free coherent sheaves $\overline{\cE}$ on $S$ parametrized by $\boldc$, so that
\bd
\overline{\cE}\big|_{\boldc\times D} \cong E\otimes\theo_\boldc.
\ed
\end{cor}

We include the proof of this fact, because the proof arises again, in conjunction with the $\cs$ action on $\M_S(E)$,  in Section \ref{preliminaries}.

\begin{proof}
The map 
\bd
m: \cs \times S \longrightarrow S
\ed
when restricted to $\cs\times (S\setminus \sigma)$ extends to a map 
\bd
\overline{m}: \boldc \times (S\setminus\sigma) \longrightarrow S\setminus\sigma.
\ed
One gets a sheaf $\tilde{\cE}:= \overline{m}^*\cE$, whose restriction to $\cs\times(S\setminus\sigma)$ is canonically isomorphic to $m^*\cE$.  One clearly gets a framing
\bd
\tilde{\cE} \xrightarrow{\overline{m}^*\phi} E\otimes\theo_\boldc
\ed
from the framing $\cE \xrightarrow{\phi} E$.  
Restricting $\tilde{\cE}$ to $\boldc\times(S\setminus(\sigma\cup D))$ and $\overline{m}^*\cE$ to $\cs\times(S\setminus D)$ gives exactly the input data for the Extension Lemma, and so determines a unique extension of $\tilde{\cE}$ to $\boldc\times(S\setminus D)$.  Finally Zariski gluing now pieces together $\overline{\cE}$ from $\tilde{\cE}$ and its extension across $\sigma$.
\end{proof}

A slightly more conceptual restatement of this corollary is the following.

\begin{cor} Suppose $(\cE,\phi)$ is an element of $\M_S(E)(\spec\boldc)$.  Then a limit 
\bd
\lim_{\lambda\rightarrow 0}\lambda\cdot(\cE,\phi)
\ed
 exists in $\M_S(E)(\spec\boldc)$.
\end{cor}

One has also a completeness property for the fixed-point set $(\M_S(E))^{\cs}$.

\begin{prop}\label{fpscompleteness}(\cite{framedlocalize})
Suppose $R$ is a DVR, $K$ is its field of fractions, and 
\bd
(\cE_K, \phi_K)\in (\M_S(E))^{\cs}(\spec K).
\ed
Then $(\cE_K,\phi_K)$ is the restriction of some family 
\bd
(\cE_R, \phi_R)\in(\M_S(E))^{\cs}(\spec R).
\ed 
\end{prop}

We end this section with a tool that allows us to prove separatedness of the stack $\M_S(E)$ in some cases.  We continue to use the notation of the earlier parts of this section without further comment.

\begin{uniqlemma}
Suppose that $\cE_R$ and $\cE_R'$ are $R$-flat families of torsion-free rank two $\cA_R$-modules and that, for the associated families of torsion-free sheaves on $\bspec \Sym B$, there is some neighborhood of $C_R$ in $\bspec\Sym B$ on restriction to which the associated families of sheaves are locally free.  Suppose furthermore that there are given isomorphisms
\begin{align}
\cE_R \otimes K\rightarrow & \cE_R'\otimes K \text{\hspace{1em}and}\\
\cE_R\otimes \widehat{\cA}_R \rightarrow & \cE_R'\otimes \widehat{\cA}_R
\end{align}
that yield a commutative diagram
\begin{equation}\label{uniqcd}
\xymatrix{ & & \cE_R\otimes K \ar[dd]\ar[dr] & \\ & & & \cE_R'\otimes K\ar[dd]\\ \cE_R\otimes \widehat{\cA}_R \ar[rr]\ar[dr] & & \cE_R\otimes\widehat{\cA}_R \otimes K \ar[dr] & \\ & \cE_R'\otimes \widehat{\cA}_R \ar[rr] & & \cE_R'\otimes\widehat{\cA}_R\otimes K\hspace{.5em}.}
\end{equation}
Then the given morphism $\cE_R\otimes K\rightarrow \cE_R'\otimes K$ induces an isomorphism $\cE_R\cong \cE_R'$ via the natural inclusions.
\end{uniqlemma}

\begin{proof}
We will show that the images of $\cE_R$ and $\cE_R'$ in $\cE_R\otimes K$ and $\cE_R'\otimes K$, respectively, are identified under the given isomorphism.

Locally along $C$ we may trivialize $B$, and then over an open set $U_R\subset C_R$ we may replace $\cA_R$, $\widehat{\cA}_R$ by $T_R[t]$, $T_R\ldb t\rdb$ respectively.  Choose $f\in T_R[t]\setminus (t)$ so that $(\cE_R)_f$ and $(\cE_R')_f$ are free modules over $T_R[t]_f$; we may choose such an $f$ (assuming we have chosen $U_R$ sufficiently small) because $\cE_R$ and $\cE_R'$ are locally free near $C_R$.

\begin{claim}\label{claiminuniq}
\mbox{}\noindent
\begin{enumerate}
\item[a.] $(\cE_R)_f$  is the intersection of $(\cE_R)_f\otimes K$ and 
\mbox{$(\cE_R)_f\underset{T_R[t]_f}{\otimes}T_R\ldb t\rdb$} in \newline\mbox{$(\cE_R)_f\underset{T_R[t]_f}
{\otimes}T_R\ldb t\rdb\underset{R}{\otimes}K$}.
\item[b.] $(\cE_R')_f$ is the intersection of $(\cE_R')_f\otimes K$ and 
\mbox{$(\cE_R')_f\underset{T_R[t]_f}{\otimes} T_R\ldb t\rdb$} in \newline\mbox{$(\cE_R')_f\underset{T_R[t]_f}{\otimes} T_R\ldb t\rdb\underset{R}{\otimes} K$}.
\end{enumerate}
\end{claim}
\begin{proof}
By the choice of $f$ it is enough to show that
\bd
T_R[t]_f =\left( T_R[t]_f\otimes K\right) \cap T_R\ldb t\rdb.
\ed
So, write an element of the intersection as
\bd
\frac{g}{\pi^jf^k} = s, \text{\hspace{.5em}where\hspace{.5em}} g\in T_R[t]\text{\hspace{.5em}and\hspace{.5em}} s\in T_R\ldb t\rdb.
\ed
Then
\bd
g=\pi^jf^ks\in \pi^j T_R\ldb t\rdb.
\ed
Therefore, all coefficients of $g$ are divisible by $\pi^j$, and $\ds\frac{g}{\pi^j}\in T_R[t]$ after all; consequently, $\ds\frac{g}{\pi^jf^k}\in T_R[t]_f$ as desired.  This proves the claim. \end{proof}

Claim \ref{claiminuniq}, together with the commutative diagram \eqref{uniqcd}, gives an isomorphism
\bd
(\cE_R)_f\rightarrow (\cE_R')_f
\ed
so that the diagram
\begin{equation}
\xymatrix{ & & \cE_R\otimes K \ar[dd]\ar[dr] & \\ & & & \cE_R'\otimes K\ar[dd]\\ (\cE_R)_f \ar[rr]\ar[dr] & & (\cE_R)_f \otimes K \ar[dr] & \\ & (\cE_R')_f \ar[rr] & & (\cE_R')_f\otimes K}
\end{equation}
commutes.  Then by Proposition 6 of \cite{MR51:510}, the inclusions
\begin{align*}
\cE_R\rightarrow & \cE_R\otimes K,\\
\cE_R'\rightarrow & \cE_R'\otimes K
\end{align*}
induce an isomorphism of the  $\cA_R$-modules $\cE_R$ and $\cE_R'$ over the open set $U_R$ of $C_R$.  The isomorphisms so obtained locally along $C_R$ arise from the restrictions of a single diagram
\begin{equation*}
\xymatrix{\cE_R \ar[d] & \cE_R'\ar[d] \\ \cE_R\otimes K \ar[r] & \cE_R'\otimes K}
\end{equation*}
over $C_R$, so they are compatible on intersections of the open sets of the form $U_R$ and consequently induce an isomorphism $\cE_R\rightarrow \cE_R'$ over all of $C_R$.
\end{proof}

\section{Representability of the Stack}\label{chrepresentability}

In this section, we prove representability and good behavior of the stack $\M_S(E)$ in some special cases.  In fact, one can show that this stack may be replaced by a reasonably well-behaved scheme with the same rational homology under very general conditions, but since we will not use this fact further in this work we omit it (the interested reader may find this result in \cite{framedlocalize}), focusing instead on the behavior of the stack itself.  Some related representability results may be found in \cite{MR95g:14013} and \cite{framedrep}.

We begin by reviewing the description of polarizations of the surface 
$S\xrightarrow{\pi} C$, since we will use this description below.  First, 
recall that $\operatorname{Pic} S\cong \operatorname{Pic} C \oplus 
\mathbf{Z}$, where the factor $\mathbf{Z}\subset \operatorname{Pic} S$ 
may be generated by a section of the projection $\pi$.  Suppose $\tau$ is
a section of $\pi$ whose self-intersection $\tau^2$ is minimal among 
self-intersections of sections of $\pi$; note that $\tau^2$ is then the
same as
\bd
-\operatorname{max} \lbrace 2\deg B - \deg L \Big| \text{$B$ is a 
subbundle of $L\oplus\theo$}\rbrace.
\ed
If $\deg L\geq 0$, then this is just $-\deg L$, and the divisor at infinity
$D$ is such a section.  In this setting, one has (see \cite{MR99c:14056}) that $H=aD+bf$ (where $f$ is the numerical class of a fiber) is ample if and only if $a>0$ and $b>\frac{a}{2}\deg L$.

\subsection{Elliptic Curves Equipped with Stable Bundles}

If $L$ has degree zero, one need not impose the usual stability or semistability conditions to get an open substack that is an honest moduli {\em space}; one has rather the following fact.

\begin{prop}\label{norissuggestion}
Suppose $C$ is a curve of any genus, $L$ has degree zero, $E$ is a rank two stable vector bundle on $D$, and $(\cE, \phi)\in \M_S(E)(\spec\boldc)$.  Then $\cE$ is a stable sheaf (in the sense of Mumford-Takemoto) for an appropriate choice of polarization $H$ of $S$ that depends only on the Chern classes of $\cE$.  
\end{prop}
\begin{proof}
Suppose $\cL\subset \cE$ is a saturated rank one subsheaf of $\cE$, with quotient $\cE/\cL = Q$.

Suppose first that the restriction of $\cE$ to a generic fiber of $S\xrightarrow{\pi} C$ is of the form $\cE\big|_f \cong \theo_{\pline}(d)\oplus\theo_{\pline}(d)$.  Write $\cL = I \otimes \theo (d'\sigma)\otimes \pi^*B$, where $d'$ is the generic fiber degree of $\cL$ (that is, the degree of the restriction of $\cL$ to a generic fiber of $\pi$), $B$ is a line bundle on $C$, and $I\subseteq \theo_S$ is an ideal of finite colength.  Then the restriction to a generic fiber $f$ of $\pi$ gives
\bd
0\rightarrow \cL\big|_f \rightarrow \cE\big|_f \rightarrow Q\big|_f \rightarrow 0,
\ed
hence an injection $\theo(d')\subset \theo(d)\oplus\theo(d)$; consequently, \mbox{$d'\leq d \leq \deg Q\big|_f$}.  

Moreover, since $\cE$ is locally free near $D$ and $\cL$ is saturated, one has
\bd
\operatorname{supp}(\theo/I)\cap D = \emptyset,
\ed
 and so one gets
\bd
B = \cL\big|_D \subset \cE\big| _D = E.
\ed
It follows that $\deg B < \deg E$ since $E$ is stable.

In particular, one has $c_1(\cL)\cdot f \leq c_1(Q)\cdot f$ and $c_1(\cL)\cdot D < c_1(Q)\cdot D$; consequently, $\cE$ is stable for any polarization $H=aD + bf$ for which $a, b>0$.  Because $\sigma$ is equivalent to $D+\pi^*c_1(L)$, we get
\bd
H = a(\sigma - (\deg L)\cdot f) + bf = a\sigma + bf
\ed
in homology.  So $\cE$ is already stable for any polarization $H= a\sigma + bf$ for which $a, b>0$.

Now, suppose $\cE\big| _f \cong \theo (d)\oplus \theo (d')$ generically, where $d>d'$.  If $\cL\subset \cE$ is destabilizing, then $\cL\big| _f \cong \theo (d)$ generically, for otherwise, as before, the restrictions of $\cL$ to both $D$ and a generic fiber $f$ are of lower degree than the restrictions of the quotient $\cE/\cL$, proving that $\cL$ is not destabilizing.

Therefore, in this case $\cL\subseteq \bar{\cL} \subset \cE$, where $\bar{\cL}$ is the canonically defined rank one subsheaf of $\cE$; thus, we may assume that $\cL=\bar{\cL}$.  Writing 
\bd
\cL = I_1\otimes \theo (d\sigma) \otimes \pi^*B_1
\ed
and
\bd
Q = \cE/{\mathcal L} = I_2\otimes \theo(d'\sigma)\otimes \pi^*B_2,
\ed
where $I_1$ and $I_2$ are ideals of finite colength and $B_1$ and $B_2$ are line bundles on $C$, we get
\bd
(c_1{\mathcal L})\cdot H = \left(d\sigma + (\deg B_1)f\right)\cdot (a\sigma + bf) = a\deg B_1 + bd,
\ed
and similarly 
\bd
(c_1Q)\cdot H = a\deg B_2 + bd'.
\ed
So for fixed $d, d'$, and $b$, by choosing $a$ sufficiently large one gets $(c_1Q)\cdot H > (c_1{\mathcal L})\cdot H$.  The lemma that follows proves that one can choose the coefficient $a$ sufficiently large relative to $b$ to work uniformly for all sheaves $\cE$ of fixed Chern classes $c_1$ and $c_2$.
\end{proof}

\begin{lemma}\label{lemmaoffinite} Fix $c_1$ and $c_2$.  There are only finitely many choices of $d$ and $d'$ for which there exists a sheaf $\cE$ on $S$ satisfying
\begin{enumerate}
\item $\cE\big| _D \cong E$, 
\item $c_1\cE =c_1$,
\item $c_2\cE = c_2$, and 
\item $\cE\big| _f \cong \theo(d)\oplus \theo(d')$ for $f$ the generic fiber of $S\overset{\pi}{\longrightarrow} C$.
\end{enumerate}
\end{lemma}
\begin{proof}
We may assume without loss of generality that $E$ has degree $0$ or $1$ and that the restriction of $\cE$ to a generic fiber of $\pi$ has degree \mbox{$F = d+d'\leq 0$}: these are determined by $c_1$, and by using twists by powers of $\theo_S(\sigma)$ and pullbacks of line bundles on $D$, we may move any component of a moduli stack $\M_S(E)$ isomorphically to a component of a moduli stack $\M_S(E')$ with the appropriate degrees. 

Suppose now that $\cE$ satisfies the four requirements in the statement of Lemma \ref{lemmaoffinite}.  Assume, without loss of generality, that 
\bd
\cE\big|_f \cong \theo (d) \oplus \theo (d')
\ed
with $d> d'$.  One gets the canonical short exact sequence 
\bd
0 \rightarrow I_1 \otimes \theo(d\sigma) \otimes \pi^*B_1 \rightarrow \cE \rightarrow I_2\otimes \theo (d'\sigma)\otimes \pi^*B_2 \rightarrow 0.
\ed
Then $\cE\big|_D \cong E$ implies that
\bd
\deg B_1 + \deg B_2  = \deg E.
\ed
Moreover,
\begin{equation}
\begin{align}
c_2\cE & = (d\sigma + f\deg B_1) \cdot (d'\sigma + f\deg B_2) + c_2 I_1 + c_2 I_2\\
& = d'\cdot\deg B_1 + d\cdot \deg B_2 + c_2I_1 + c_2 I_2\\
& = d'\cdot\deg B_1 + d(\deg E - \deg B_1) + c_2I_1 + c_2 I_2\\
& = d\cdot \deg E + (d'-d)\deg B_1 + c_2 I_1 + c_2 I_2\\
& = d\cdot \deg E + (F-2d)\deg B_1 + c_2 I_1 + c_2 I_2.\label{lastone}
\end{align}
\end{equation}
If $\deg E =0$ then $\deg B_1 <0$; combining this with our assumption that $F-2d <0$, we see that all terms of Equation \eqref{lastone} are nonnegative and that the second term increases linearly with $d$, which gives a bound on $d$ since $c_2\cE$ is fixed.  Similarly, if $\deg E = 1$ then $\deg B_1 \leq 0$, and consequently all terms of Equation \eqref{lastone} are nonnegative with the first term increasing linearly with $d$; thus we again obtain a bound on $d$.  This proves the lemma.
\end{proof}

Proposition \ref{norissuggestion} has the following pleasant consequence when the curve $C$ has genus 1.

\begin{cor}\label{noricor}
If $C$ is an elliptic curve, $E$ is stable, and $\deg L =0 $, then $\M_S(E)$ is represented by a smooth quasiprojective variety.
\end{cor}

\begin{proof}
Suppose that $(\cE, \phi)$ is an $E$-framed pair and $\cE$ is stable.  We
will prove that then the pair $(\cE,\phi)$ is a stable pair in the
sense of Huybrechts--Lehn (\cite{MR95i:14015}, \cite{MR96a:14017}) for the ample divisor $H$ chosen in the lemma above and the auxiliary datum of the polynomial $\delta$; because our stack is then an open substack of their moduli stack of stable pairs and they prove the representability of their stack by a quasiprojective variety, the representability of $\M_S(E)$ follows.  

In what follows, we choose $\delta := p_{E}$, the Hilbert polynomial of $E$.  Then the Hilbert polynomial of the pair $(\cE,\phi)$ as defined by Huybrechts--Lehn is exactly
\begin{align*}
p_{(\cE,\phi)} & := p_{\cE} - \delta\\
 & = p_{\cE (-D)}.
\end{align*}

Suppose, then, that $\cL\subset \cE$ is a rank one subsheaf.  If $\cL\subset \cE(-D)$, then the Hilbert polynomial associated to the pair $(\cL, \phi\big| _{\cL})$ is
\bd
p_{(\cL,\phi\big| _{\cL})} := p_{\cL},
\ed
the Hilbert polynomial of $\cL$.  Since $\cE(-D)$ is stable, we have
\bd
p_{(\cL,\phi\big| _{\cL})} = p_{\cL} < \frac{1}{2}p_{\cE(-D)} = \frac{1}{2}p_{(\cE,\phi)}
\ed
as desired.  If, on the other hand, $\cL$ is not contained in $\cE(-D)$, then
\bd
p_{(\cL,\phi\big|_{\cL})} = p_{\cL} - \delta,
\ed
and
%% \begin{align*}
\bd
p_{(\cL,\phi\big|_{\cL})} = \bigl( p_{\cL} - \delta\bigr) < \biggl( \frac{1}{2}p_{\cE} -\delta\biggr)
 \leq \biggl( \frac{1}{2}p_{\cE} - \frac{1}{2}\delta\biggr) = \frac{1}{2}p_{(\cE,\phi)}
\ed
%% \end{align*}
as desired.

Huybrechts--Lehn prove in \cite{MR96a:14017} that $M_S(E)$ is smooth whenever the hyper-Ext obstruction group $\bExt^2_S (\cE, \cE\overset{\phi}{\longrightarrow} E)$ is zero for all pairs $(\cE,\phi)$ in $M_S(E)$.  Since the complex $\cE\overset{\phi}{\longrightarrow} E$ is quasi-isomorphic to $\cE(-D)$ in our setting, one obtains
\bd
\bExt^2_S (\cE,\cE\overset{\phi}{\longrightarrow}E) \cong \Ext^2_S (\cE,\cE(-D)).
\ed
Now by Serre duality, 
\begin{equation}
\begin{align}
\Ext^2_S(\cE,\cE(-D))^* &\cong \Ext^0_S (\cE (-D), \cE\otimes K_S)\\
 & = \Ext^0_S(\cE,\cE(D)\otimes\theo(-\sigma-D)\otimes\pi^*K_C)\\
& = \Ext^0_S (\cE,\cE(-\sigma))\\
& = 0 \text{\hspace{1em} because $\cE$ is simple.}
\end{align}
\end{equation}
This completes the proof of the corollary.
\end{proof}

\subsection{Elliptic Curves Equipped with Polystable Bundles}

There is a family of special cases in which framing by a polystable bundle gives a reasonable moduli space, on which it will be straightforward to compute the rational homology.

\begin{thm}\label{polystableisgood}
Suppose that $C$ is an elliptic curve, that $L_1$ and $L_2$ are line bundles on $C$ both of the same arbitrary degree, and that $L$ is a line bundle of degree zero on $C$.  Suppose further that 
\begin{enumerate}
\item $L^k$ is nontrivial for all $k\in\mathbf{Z}\setminus \lbrace 0 \rbrace$; and
\item $L^k\otimes (L_1\inv \otimes L_2)$ is nontrivial for all $k\in\mathbf{Z}$.
\end{enumerate}
Then if $S=\mathbf{P} (L\oplus \theo)$, the moduli stack $\psstack$ is a smooth separated scheme.
\end{thm}

\begin{proof}
We may assume that $\deg L_1 = \deg L_2 = 0$, since if we begin with $\deg L_1 = \deg L_2 = d$ then for any point $p\in C$, $L_1' = L_1(-dp)$ and $L_2' =L_2(-dp)$ will satisfy the assumptions.

We begin by proving that, for any $\spec\boldc$-valued point $(\cE,\phi)$ of $\psstack$, we have
\bd
\Ext^0_S (\cE,\cE(-D))= \Ext^2_S(\cE,\cE (-D)) = 0.
\ed
In this computation, we will use repeatedly that
\begin{align*}
\Ext^2_S(\cE,\cE(-D))^* &\cong \Ext^0_S(\cE(-D),\cE\otimes\theo(-\sigma -D))\\
&= \Ext^0_S (\cE,\cE(-\sigma)).
\end{align*}

Suppose first that, if $f$ denotes the generic fiber of the projection $\pi :S\rightarrow C$, we have
\bd
\cE\big| _f \cong \theo(l)\oplus\theo(l).
\ed
Then any endomorphism of $\cE\big|_f$ that vanishes at infinity is zero; hence any $e\in\Ext^0_S(\cE,\cE(-D))$ is zero on the generic fiber of $\pi$ and consequently is zero on all of $S$.  Similarly, $\Ext^0_S(\cE,\cE(-\sigma))=0$ and hence $\Ext^2_S(\cE,\cE(-D))=0$.

Suppose next that 
\bd
\cE\big|_f \cong \theo (l)\oplus\theo(l'\sigma)
\ed
with $l > l'$.  Then, as in Section \ref{exactseqtypeU}, we get a canonical exact sequence
\bd
0\rightarrow \pi^*B\otimes\theo(l\sigma)\otimes I_1 \rightarrow \cE \rightarrow \pi^*(B\inv\otimes L_1\otimes L_2)\otimes\theo (l'\sigma)\otimes I_2\rightarrow 0.
\ed
Now
\begin{align*}
 \Ext^0 \Bigl(\pi^*(B\inv\otimes L_1  \otimes L_2)\otimes \theo & (l'\sigma) \otimes I_2, \pi^*B\otimes\theo (l\sigma)\otimes I_1\Bigr)\\
& \subseteq \Ext^0 \left(\pi^*(B\inv\otimes L_1\otimes L_2)\otimes\theo (l'\sigma), \pi^*B\otimes\theo(l\sigma)\right)\\
& = H^0\left(S,\pi^*B^2\otimes L_1\inv\otimes L_2\inv \otimes\theo((l-l')\sigma)\right) \\
& = H^0 \left(C, B^2\otimes L_1\inv\otimes L_2\inv\otimes \sum_{j=0}^{l-l'}L^j\right).
\end{align*}

If $\deg B <0$, this group is zero, while if $\deg B = 0$, then, since $L_1\oplus L_2$ is semistable of degree zero, we have either $B=L_1$ or $B=L_2$, and in either case the second assumption of the theorem implies that this group is zero.

Now if the canonical exact sequence is nonsplit, the above vanishing implies that $\cE$ is simple, hence $\Ext^0_S(\cE,\cE(-D)) = 0$ as desired.  If, on the other hand, the canonical exact sequence splits, then $B$ and $L_1\otimes L_2\otimes B\inv$ are line subbundles of the polystable bundle $L_1\oplus L_2$, and so $\deg B = 0$ and either $B=L_1$ or $B=L_2$.  This implies the $H^0$ vanishing just as before.  To prove that $\Ext^0_S(\cE,\cE(-D))=0$ in this case, we must compute also
\begin{align*}
\Ext^0 \Bigl( \pi^*B\otimes\theo (l\sigma) \otimes I_1, \pi^* & (B\inv\otimes L_1\otimes L_2)\otimes \theo (l'\sigma)\otimes I_2\Bigr)\\
& \subseteq \Ext^0 \left( \pi^*B \otimes\theo (l\sigma),\pi^*(B\inv \otimes L_1\otimes L_2)\otimes\theo (l'\sigma)\right)\\
& = H^0\left( S, \pi^*(B^{-2}\otimes L_1\otimes L_2) \otimes \pi_*\theo ((l'-l)\sigma)\right)\\
& = 0
\end{align*}
since $l'<l$.  Consequently,

\begin{multline}
\Ext^0(\cE,\cE) = \End \bigl(\pi^*B\otimes\theo(l\sigma)\otimes I_1\bigr)\\
\bigoplus \End\left( \pi^* (B\inv \otimes L_1\otimes L_2)\otimes \theo (l'\sigma)\otimes I_2\right)
\end{multline}
in this case.  But the two sheaves in this formula are of rank one, hence
\bd
\Ext^0(\cE,\cE) = \boldc^2
\ed
and the restriction map along $\sigma$ or along $D$ is just the identity of $\boldc^2$.  As a result 
\bd
\Ext^0_S(\cE,\cE(-D)) = \Ext^2_S(\cE,\cE(-D))=0
\ed
here as well.

As a consequence of the vanishing of $\Ext^0_S(\cE,\cE(-D))$ for every $\spec\boldc$-point $(\cE,\phi)$ of $\psstack$, we find that every object of $\psstack$ over $\spec\boldc$ is rigid.  But then by semicontinuity every object of $\psstack$ is rigid as well, and by Propositions 4.4 and 1.4.1.1 of \cite{LMB}, the stack $\psstack$ is in fact represented by a locally finitely presented algebraic space.

The deformation theory arguments used in \cite{MR96a:14017} to obtain a smoothness criterion for moduli spaces of stable framed sheaves are arguments purely about the smoothness of the moduli functor, and hence carry over {\em mutatis mutandis} to our setting; in particular, if
\bd
\bExt^2_S (\cE,\cE\xrightarrow{\phi} L_1\oplus L_2) = 0
\ed
for all pairs $(\cE,\phi)$, then $\psstack$ is smooth.  But since $\phi$ induces an isomorphism
\bd
\cE\big| _D \cong L_1\oplus L_2,
\ed
we have that the complex $\cE\xrightarrow{\phi}L_1\oplus L_2$ is quasi-isomorphic to $\cE(-D)$ concentrated in degree zero, and so
\bd
\bExt^2_S (\cE,\cE\xrightarrow{\phi} L_1\oplus L_2) = \Ext^2_S (\cE,\cE (-D)) = 0
\ed
for all pairs $(\cE,\phi)$.  As a result, $\psstack$ is smooth.

Next, we need to prove that $\psstack$ is separated; so, writing $E=L_1\oplus L_2$, suppose that $(\cE_K,\phi_K)$ is a family of $E$-framed torsion-free sheaves parametrized by $\spec K$, where $K$ is the field of fractions of a DVR $R$.  We proceed in several steps.

\begin{claim}
Suppose $A$ is an affine scheme and $(\cE_A,\phi_A)$ is an $E$-framed family parametrized by $A$. 
Then $\cE_A\big|_{\widehat{D}_A}$ is isomorphic to $\theo_{\widehat{D}_A}\underset{\theo_{D_A}}{\otimes} E_A$; furthermore, this isomorphism is completely determined by $\phi_A$.
\end{claim}
\begin{proof}
Suppose we have shown that $\cE_A\big|_{D^{(n)}_A}$ is isomorphic to $\theo_{D^{(n)}_A}\otimes E_A$.  By Proposition 1.4 of \cite{MR34:6796}, the obstruction to uniqueness of an extension of $\theo_{D^{(n)}_A}\otimes E_A$ to a bundle over $D_A^{(n+1)}$ is a class in 
\begin{align*}
H^1\Bigl(\End E_A & \otimes ( I^n_D/I^{n+1}_D)\Bigr) \cong H^1\left(\End E\otimes (I_D^n/I_D^{n+1})\right)\otimes A\\
& =\left[ H^1(L^n) \oplus H^1(L^n)
\oplus H^1(L_1\inv\otimes L_2\otimes L^n)\oplus H^1(L_2\inv\otimes L_1\otimes L^n)\right]\otimes A,
\end{align*}
which vanishes whenever $n\geq 1$ by the assumptions of the theorem.  So 
\bd
\cE_A\big|_{D_A^{(n+1)}} \cong \theo_{D_A^{(n+1)}}\otimes E_A.
\ed
  Moreover, since
\begin{align*}
\Hom \Bigl(E_A, E_A\otimes (I_D^n/I_D^{n+1})\Bigr) & = \Bigl[ H^0(L^n)\oplus H^0(L^n)\oplus H^0(L_1\inv\otimes L_2\otimes L^n)\\
&  \hspace{5em}\oplus H^0(L_2\inv\otimes L_1\otimes L^n)\Bigr]\otimes A\\
& =0
\end{align*}
whenever $n\geq 1$, we see that the isomorphism
\bd
\cE_A\big| _{D_A^{(n+1)}} \cong \theo _{D_A^{(n+1)}}\otimes E_A
\ed
is uniquely determined by $\phi_A$.
\end{proof}

The same computation as in the proof of this claim gives
\begin{align}\label{num1}
\End(\theo_{\widehat{D}_R}\otimes E_R) & \cong \End(E_R) \cong \End(E)\otimes R, \text{\hspace{2em}and}\\ \label{num2}
\End(K\otimes\theo_{\widehat{D}_R}\otimes E_K) & \cong \End(E_K) \cong \End(E)\otimes K.
\end{align}
Suppose $(\cE_R,\phi_R)$ and $(\cE_R',\phi_R')$ are $E$-framed families parametrized by $\spec R$ that are equipped with isomorphisms
\begin{align}
(\cE_R,\phi_R)\big|_{S_K} & \xrightarrow{\psi} (\cE_K,\phi_K) \text{\hspace{2em}and}\\
(\cE_R',\phi_R') \big|_{S_K} & \xrightarrow{\psi'} (\cE_K,\phi_K).
\end{align}
Restricting to $\widehat{D}_R$ and identifying each of $\cE_R\big|_{\widehat{D}_R}$, $\cE_R'\big|_{\widehat{D}_R}$ uniquely with $\theo_{\widehat{D}_R}\otimes E_R$, we obtain an automorphism
\bd
a = \left( (\psi')\inv\circ \psi\right) \big|_{\widehat{D}_R\otimes K}
\ed
of $\bigl( K\otimes\theo_{\widehat{D}_R}\bigr)\otimes E_K$ whose restriction to $E_K$ is the identity.

\begin{claim}
The automorphism
$a$ is the restriction of the identity automorphism of $\theo_{\widehat{D}_R}\otimes E_R$.
\end{claim}
\begin{proof}
This follows immediately from the description in Equations \eqref{num1} and 
\eqref{num2}.\nolinebreak\end{proof}
Now, restrict to $S_R\setminus \sigma_R$.  We then have sheaves of $\cA_R$-modules $\cE_R$, 
$\cE_R'$ where 

\noindent
$\cA = \Sym L$, together with embeddings
\begin{align}
\cE_R \rightarrow & \theo_{\widehat{D}_R}\otimes \cE_R,\\
\cE_R' \rightarrow & \theo_{\widehat{D}_R}\otimes \cE_R';
\end{align}
moreover, the isomorphisms 
\begin{align}
\theo_{\widehat{D}_R}\otimes \cE_R \rightarrow & \theo_{\widehat{D}_R}\otimes E_R \text{\hspace{2em}and}\\
\theo_{\widehat{D}_R}\otimes \cE_R' \rightarrow & \theo_{\widehat{D}_R}\otimes E_R
\end{align}
identify $\cE_R$ and $\cE_R'$ with subsheaves of $\theo_{\widehat{D}_R}\otimes E_R$, and in addition the restrictions of their images to $D_K$ coincide.  By the Uniqueness Lemma, then, the inclusions
\begin{align}
\cE_R\rightarrow & \cE_K \text{\hspace{2em}and}\\
\cE_R'\rightarrow & \cE_K
\end{align}
have the same image and thus induce an isomorphism of $\cE_R$ and $\cE_R'$ over $S_R\setminus \sigma_R$.  Because this is compatible with the isomorphism given over $S_K$, by the \mbox{Extension Lemma} we obtain an isomorphism
\bd
(\cE_R,\phi_R)\cong (\cE_R',\phi_R').
\ed
This completes the proof that $\psstack$ is separated.

We now indicate why $\psstack$ is in fact a scheme.  Suppose $(\cE,\phi)$ is a $\spec\boldc$-valued point of $\psstack$.  Let $H$ denote a polarization of $S$, and, given an $\theo_S$-module $\cF$, write
\bd
\cF(n) := \cF\otimes H^n.
\ed
Choose $n$ sufficiently large that $\cE(n)$ is globally generated and
\bd
H^1\bigl(\cE(n)\bigr) = H^2\bigl(\cE(n)\bigr) = 0;
\ed
then the same will hold for all $\cF$ occurring in pairs $(\cF,\psi)$ that lie in an open neighborhood of $(\cE,\phi)$ in $\psstack$. 

We saw above that $\cE\big|_{D^{(m)}} \cong \theo_{D^{(m)}}\otimes E$ for any $m\geq 1$ and that the isomorphism is uniquely determined by $\phi$.  For $m$ sufficiently large, this isomorphism determines an injection
\bd
H^0\bigl(\cE(n)\bigr)\rightarrow H^0\bigl(\theo_{D^{(m)}}\otimes E(n)\bigr);
\ed
fix such an $m$.

We next construct a scheme representing a neighborhood of $(\cE,\phi)$ in \mbox{$\psstack$}.  There is a Grassmannian $G$ parametrizing subspaces of \mbox{$H^0\bigl(\theo_{D^{(m)}}\otimes E(n)\bigr)$} that are of {\nolinebreak[4] dimension} equal to the dimension of $H^0\bigl(\cE(n)\bigr)$.  Over $S\times G$, in addition, there is a
canonical subbundle $V$ of \mbox{$\theo_{S\times G}\otimes H^0\bigl(\theo_{D^{(m)}}\otimes E(n)\bigr)$}; from this family $V$ on $S\times G/G$, one obtains a relative Quot-scheme $\cQ$ together with a universal subsheaf $\cU$ of $p_{S\times G}^* V$ on $S\times G\times\cQ$: this Quot-scheme parametrizes subsheaves of the bundles $V_g$ ($g\in G$) on $S$ with Hilbert polynomials coinciding with that of
\bd
\operatorname{ker} \Bigl( \theo_S\otimes H^0\bigl(\cE(n)\bigr)\rightarrow \cE(n)\Bigr).
\ed

We have a surjective map
\bd
\theo_S\otimes H^0\bigl(\theo_{D^{(m)}}\otimes E(n)\bigr) \rightarrow \theo_{D^{(m)}}\otimes E(n),
\ed
and there is a locally closed subvariety $\cQ'$ of $\cQ$ parametrizing those \mbox{$\cU_{g\times q}\subset V_g$} on $S$ for which
\begin{enumerate}
\item $\operatorname{ker} \Bigl( V_g\rightarrow E(n)\Bigr)$ is contained in $\cU_{g\times q}$,
\item $V_g/\cU_{g\times q}$ is torsion-free, and
\item the induced map
\bd
\left(V_g/\cU_{g\times q}\right) \big| _D \rightarrow E(n)
\ed
is an isomorphism.
\end{enumerate}
Over $\cQ'$ we have a universal diagram
\bd
\xymatrix{V\ar[d]\ar[r] & \theo_{S\times G\times\cQ'}\otimes H^0\bigl(\theo_{D^{(m)}}\otimes E(n)\bigr) \ar[d]\\
V/\cU \ar[r]^{\psi} & E(n)_{G\times\cQ'}\hspace{.2em}.}
\ed

\begin{claim} The map \mbox{$\cQ'\rightarrow \psstack$} that is induced by the family \newline
\mbox{$\left(V/\cU (-n), \psi\big|_D (-n)\right)$} represents a neighborhood of $(\cE,\phi)$ in $\psstack$.
\end{claim}

\begin{proof}
Suppose $(\cF_U,\xi_U)$ is a $U$-flat family in $\psstack$ for which
\begin{enumerate}
\item the Hilbert polynomial of each $\cF_{U,u}$ ($u\in U$) coincides with that of $\cE$,
\item each $\cF_{U,u}(n)$ is globally generated,
\item each $H^1\bigl(\cF_{U,u}(n)\bigr) = H^2\bigl(\cF_{U,u}(n)\bigr) =0$, and
\item each $H^0\bigl(\cF_{U,u}(n)\bigr)$ injects into $H^0\bigl(\theo_{D^{(m)}}\otimes E(n)\bigr)$ under the canonical homomorphism.
\end{enumerate}
Then the family of diagrams
\bd
\xymatrix{(p_U)^*(p_U)_* \cF_U(n) \ar[d]\ar[r] & \theo_{S\times U}\otimes H^0\bigl(\theo_{D^{(m)}}\otimes E(n)\bigr)\ar[d] \\
\cF_U(n)\ar[r] & \theo_{D^{(m)}}\otimes E(n)_U}
\ed
determines a map $U\rightarrow \cQ'$.  This functor gives exactly the inverse of the map 
\mbox{$\cQ'\rightarrow \psstack$}, proving the claim.\end{proof}

This proves that $\psstack$ is locally representable by schemes, and so is itself a scheme.

In fact, $\cQ'$ is a quasiprojective variety and its image in $\psstack$ is $\cs$-invariant; hence, in particular (see Section 4 of \cite{MR51:3186}), $\psstack$ may be covered by $\cs$-invariant quasi-affine open subschemes, and thus, by \cite{MR51:3186}, $\psstack$ admits a Bia{\l}ynicki-Birula decomposition.
\end{proof}

We will use this result in Section \ref{chhomologybasis} when we compute a basis for the rational homology of $\psstack$.

\section{A Rational Homology Basis}\label{chhomologybasis}

In this section, we use localization to compute a rational homology basis for the moduli space $\M_S(L_1\oplus L_2)$ in the case in which the line bundles $L$, $L_1$, and $L_2$ are of degree zero and satisfy the assumptions of Theorem \ref{polystableisgood}.  

\subsection{The Localization Formula}

The localization theory that we use was originally developed by
Frankel (\cite{MR24:A1730}), Bia{\l}ynicki-Birula
(\cite{MR51:3186}, \cite{MR53:13228}, \cite{MR56:12020}), and
Carrell-Sommese (\cite{MR80m:32032}); see also the work of Nakajima (\cite{MR1711344}),
  Carrell-Goresky (\cite{MR85d:32063}) and Kirwan (\cite{MR89c:14072}) for useful variations of the theory.   
The theorem that we want to use here is the following: if $M$ is a nonsingular algebraic variety with $\cs$ action, one may compute the homology of $M$ in terms of that of the fixed-point set of the $\cs$ action.  More precisely, write
\bd
M^{\cs} = \coprod_{\gamma} F_{\gamma},
\ed
a disjoint union of connected components, and let
\bd
S_{\gamma} = \lbrace m\in M \arrowvert \lim_{\lambda \rightarrow 0} \lambda\cdot m \in F_{\gamma}\rbrace.
\ed
Let $m_{\gamma}$ denote the complex codimension of $S_{\gamma}$ in $M$ (this is the same as the dimension of the negative weight space for the $\cs$ action in a fiber of the normal bundle to $F_{\gamma}$ in $M$).  Then
\bd
H_i(M;\mathbf{Q}) \cong \bigoplus_{\gamma} H_{i-2m_{\gamma}}(F_{\gamma};\mathbf{Q}).
\ed
Originally this formula was proven assuming that $M$ is a projective
variety or compact K\"{a}hler manifold.  However, it has since been
observed (see Chapter 5 of \cite{MR1711344}) that the same proof applies provided that $M$ is a smooth, Hausdorff K\"{a}hler manifold for which the flow lines 
\bd
\lbrace \lambda\cdot m \arrowvert \lambda\in (0,1]\subset \cs\rbrace
\ed
are contained in compact subsets of $M$ for all $m\in M$, using a Morse-theoretic decomposition of the manifold.  In fact, the same proof also applies when one uses the (algebraically defined) Bia{\l}ynicki-Birula decomposition on a smooth separated scheme that is covered by $\cs$-invariant quasiaffine open subschemes, provided a certain additional condition on the strata is satisfied; in Section \ref{preliminaries}, we describe this condition and prove that it is indeed satisfied by $\psstack$.  This is the setting in which we apply the localization formula to compute rational homology.

\subsection{Rational Homology of $\psstack$}

We will use the same technique used by Nakajima in \cite{MR95g:58031}: the bundle $L_1\oplus L_2$ has an action by $\cs\times\cs$ obtained by scaling in the factors, and this induces an additional $\cs\times\cs$ action on $\psstack$ by scaling the framing.  By Proposition \ref{fpscompleteness}, the components of the fixed-point set are complete, and hence we may use localization for this action to reduce the computation to computation of the homology of the fixed-point set for the entire $\cs\times\cs\times\cs$ action.  Moreover, the coefficient giving the shift in homology grading, that is, the rank of the negative normal bundle to a component of the fixed-point set, will be the same whether one computes as though localization were iterated (that is, first localizing along the $\cs$ fixed-point set and then localizing for the additional $\cs\times\cs$ action inside here) or performed for $\cs\times\cs\times\cs$ all at once.  Thus, we compute formally as though we were applying 
localization to the $(\cs)^3$ action, although we will never consider any sort of 
completeness result for the entire $(\cs)^3$ action.

The decomposition of rational homology that we will obtain is the following.  Let $\psstack (c_1,c_2)$ denote the substack of $\psstack$ consisting of pairs $(\cE,\phi)$ for which $c_1\cE = c_1$ and $c_2\cE = c_2$.  Then
\bd
\psstack = \coprod_{c_1,c_2} \psstack (c_1,c_2).
\ed

\begin{thm}\label{homologybasisthm}
Fixing 
\bd
c_1 = -l'\sigma, \text{\hspace{2em}} c_2\geq 0,
\ed
one has
\begin{multline}
H_k\Bigl(\psstack (c_1,c_2);\mathbf{Q}\Bigr)\cong \\
 \bigoplus_{l\in\mathbf{Z}}\bigoplus_{\substack{\alpha,\beta\\  |\alpha | + |\beta | = c_2}} H_{k-2d(\alpha, \beta, l, l')}\left(\Symp^\alpha C\times \Symp^\beta C ;\mathbf{Q}\right).
\end{multline}
Here $d(\alpha, \beta, l, l')$, which depends on partitions
\bd
\alpha = (1^{a_1}2^{a_2}3^{a_3}\dots)\text{\hspace{1em}and\hspace{1em}}
\beta = (1^{b_1}2^{b_2}3^{b_3}\dots),
\ed
is given by 
\begin{equation}\label{shiftindex}\begin{split}
d(\alpha, \beta, l, l') = \biggl[ |\alpha | - & \ell (\alpha) +\sum_{j>i\geq 0} a_j(1 - \delta_{l'+2l-i-1,0})  \biggr]\\
+ & \biggl[ |\beta | - \ell (\beta) + \sum_{j>i\geq 0}b_j (1 - \delta_{l'+2l+i,0}) \biggr],
\end{split}\end{equation}
where $\ds|\alpha |  = \sum_{i\geq 1} ia_i$ and $\ds\ell (\alpha) = \sum_{i\geq 1} a_i$.
Moreover,
\begin{equation}\label{decompositionformula}\begin{split}
H_* \Bigl(\Symp^{\alpha}C\times\Symp^{\beta}C;\mathbf{Q}\Bigr) \cong
   \Biggl[ \bigotimes_{a_i\neq 0} & \Bigl( H_*(\boldp ^{a_i-1};\boldq)\otimes H_*(C;\boldq)\Bigr)\Biggr] \\
 &\bigotimes \Biggl[ \bigotimes_{b_j\neq 0}\Bigl( H_*(\boldp ^{b_j-1};\boldq)\otimes H_*(C;\boldq)\Bigr)\Biggr].
\end{split}\end{equation}
\end{thm}

\subsection{Preliminaries}\label{preliminaries}

We begin by clarifying the nature of the fixed-point set and then confirming that the localization formula does indeed apply to $\psstack$.

\begin{lemma}
Let $E$ denote $L_1\oplus L_2$.
Suppose $(\cE,\phi)\in\M_S(E)(\spec\boldc)$ satisfies $m_\lambda ^* (\cE,\phi) \cong (\cE,\phi)$ for all $\lambda\in\cs$.  Then $\cE$ admits a $\cs$-equivariant structure that makes $\phi$ $\cs$-equivariant for the trivial action on $L_1\oplus L_2$.
\end{lemma}
\begin{proof}
Let 
\bd
p:\cs\times S\rightarrow S
\ed
denote the projection to $S$, 
\bd
q:\cs\times S\rightarrow \cs
\ed
denote the projection to $\cs$ and, by abuse of notation, also the projection
\bd
q:\cs\times D \rightarrow \cs
\ed
to $\cs$.  We have
\bd
m_\lambda ^* \cE \cong \cE \text{\hspace{1em}for each $\lambda\in\cs$}
\ed
by assumption, and so $\Ext_q^0(m^*\cE,p^*\cE)$ is a vector bundle on $\cs$ with fiber over $\lambda\in\cs$ equal to $\Hom(m_\lambda ^*\cE,\cE)$.

There is a natural restriction map
\bd
\Ext_q^0 (m^*\cE,p^*\cE)\rightarrow \Ext_q^0 (m^*\cE\big| _D, p^*\cE\big| _D).
\ed
Now 
\begin{align*}
\Ext_q^0 (m^*\cE\big|_D,p^*\cE\big|_D) & \cong \Ext_q^0 (E_{\cs},E_{\cs})\\
& \cong \theo_{\cs} \otimes \Hom(E,E)
\end{align*}
via $\phi$, and, by assumption, the map
\bd
\Ext_q^0(m^*\cE,p^*\cE) \rightarrow \theo_{\cs}\otimes \Hom(E,E)
\ed
has in its image the identity element $1\otimes\operatorname{id}$.

Furthermore, this map is injective: fiberwise we need only check the following. 

\begin{claim}\label{injectivity}
$\Hom(\cE,\cE)\rightarrow \Hom(\cE\big|_D,\cE\big|_D)$ is injective.
\end{claim}

But this follows from the conjunction of
\begin{enumerate}
\item $\Hom(\cE,\cE)\rightarrow \Hom(\cE,\cE\big|_{D^{(n)}})$ is injective for $n\gg 0$, since $\cE$ is torsion-free; and
\item $\Hom\left(E,E\otimes(I_D^n/I_D^{n+1})\right) = 0$ for $n\geq 1$, and consequently 
\bd
\Hom(\cE,\cE\big|_{D^{(n)}})\cong \Hom(\cE,\cE\big|_D)\text{\hspace{1em}for all $n\geq 1$.}
\ed
\end{enumerate}

So the inverse image of $1\otimes\operatorname{id}$ in $\Ext_q^0(m^*\cE,p^*\cE)$ is uniquely determined.  As in the proof of Lemma \ref{uniqueequivariance}, uniqueness also guarantees commutativity of the diagram required for the map
\bd
m^*\cE\rightarrow p^*\cE
\ed
to define a $\cs$-equivariant structure.\end{proof}

\begin{notation}
We will continue to write $E$ in place of $L_1\oplus L_2$ and $\M_S(E)$ in place of $\M_S(L_1\oplus L_2)$ when it is convenient.
\end{notation}

\begin{lemma}
Suppose $(\cE,\phi)\in\M_S(E)(\spec\boldc)$ is invariant under the action of \mbox{$\cs\times\cs$} by scaling in the factors $L_1$ and $L_2$.  Then $\cE$ splits as
\bd
\cE=\cE_1\oplus \cE_2
\ed
with $\cE_1\big|_D\cong L_1$ and $\cE_2\big|_D\cong L_2$ via $\phi$.
\end{lemma}
\begin{proof}
By Claim \ref{injectivity}, the map
\bd
\Hom(\cE,\cE)\rightarrow \Hom(\cE,\cE\big|_D)
\ed
is injective, so, because
\bd
\End (L_1\oplus L_2) \cong \boldc^2
\ed
and the image of $\Hom(\cE,\cE)$ in this group contains every $(\lambda_1,\lambda_2)\in\cs\times\cs$, we must have
\bd
\Hom(\cE,\cE)\cong \Hom(\cE,\cE\big|_D).
\ed
Consider the elements $e_1$ and $e_2$ in $\Hom(\cE,\cE)$ that are the inverse images of the elements $(1,0)$ and $(0,1)$, respectively, of $\End(L_1\oplus L_2)$.  The elements $e_1$ and $e_2$ are idempotents whose images are subsheaves (necessarily torsion-free) of $\cE$; since their composites \mbox{$e_1\circ e_2$} and \mbox{$e_2\circ e_1$} are zero, the image of neither can have rank greater than 1, but also $e_1 + e_2 = \operatorname{id}$ and so $e_1(\cE)$, $e_2(\cE)$ must be rank one torsion-free subsheaves of $\cE$ for which 
\bd
\cE = e_1(\cE) + e_2(\cE).
\ed
 Letting $\cE_i = e_i(\cE)$ gives the desired splitting. \end{proof}

Suppose $(\cE,\phi)$ is a fixed point for $\cs\times\cs\times\cs$.  Then $\cE$ splits as
\begin{equation}\label{fpeq}
\cE = \bigl(\pi^*L_1\otimes \theo (l\sigma)\otimes I_1\bigr) \bigoplus \bigl( \pi^*L_2\otimes \theo ((-l'-l)\sigma)\otimes I_2\bigr),
\end{equation}
where $I_1$ and $I_2$ are $\cs$-invariant finite-colength ideals of $\theo_S$.

Observe that
\bd
c_1(\cE) = -l'\sigma
\ed
and
\begin{align*}
c_2(\cE) & = \bigl( \deg L_1\cdot f + l\sigma\bigr)\cdot \bigl( \deg L_2\cdot f - (l'+l)\sigma\bigr) + c_2I_1 + c_2I_2\\
 & = c_2I_1 + c_2I_2.
\end{align*}

\begin{remark}
In particular, neither $c_1\cE$ nor $c_2\cE$ depends on $l$, and hence even the part of $\psstack$ with fixed $c_1$ and $c_2$ may have rational homology groups of infinite rank (in fact, the formula of
 Theorem \ref{homologybasisthm} makes it clear that it does).  However, because our ultimate goal is the computation of the action of Hecke operators, which will only move sheaves a bounded amount with respect to the invariant $l$, and since, as we will see, the homology basis we produce is finite in each value of the invariant $l$ (for fixed $c_1$ and $c_2$), the homology basis will still allow reasonable computations for Hecke operators.
\end{remark}

We now describe the additional condition that $\psstack$ must satisfy in order for the localization formula to apply to it, and confirm that the condition is indeed satisfied.  

In the K\"{a}hler case, the components $F_{\gamma}$ of the fixed-point set of $\cs$ on the manifold $M$ are the critical submanifolds for a Morse-Bott function.  This allows one to define a partial order on the index set $\{\gamma\}$ so that
\bd
\overline{S_\gamma} \subseteq \bigcup_{\delta\geq\gamma} S_\delta;
\ed
in particular, their partial order on the index set is the one generated by the relations: $\gamma \leq \delta$ if $\overline{S_\gamma} \cap S_\delta \neq \emptyset$.
Atiyah and Bott (\cite{MR85k:14006}, pp. 536--537) describe how to use this ``stratification with appropriate ordering'' to obtain the localization formula.

\begin{remark}
The proof of the localization formula in the algebraic setting given by Bia{\l}ynicki-Birula does not require this partial ordering on the strata, but as his proof uses the Weil conjectures it is not so obvious to the author how one might try to modify his proof to apply to our noncompact schemes.
\end{remark}

The crucial property for the partial ordering used by Atiyah--Bott is the following: suppose 
$\Gamma$ is a subset of the index set for the strata that satisfies
\bd
\text{if $\gamma\in\Gamma$ and $\mu\leq\gamma$ then $\mu\in\Gamma$,}
\ed
for which $\ds\bigcup_{\gamma\in\Gamma}S_\gamma$ is the corresponding open submanifold of $M$. If
 $\delta$ is an index that is minimal among those indices that do not lie in $\Gamma$, then 
$S_\delta$ is a closed 
submanifold of
\bd
\bigcup_{\gamma\in\Gamma} S_\gamma \cup S_\delta.
\ed

Our scheme $\psstack$ is the union of the increasing sequence $\big\{\psstack_m\big\}$ of open subschemes that parametrize pairs $(\cE,\phi)$ for which
\begin{equation}\label{splittype}
\cE\big|_f \cong \theo(l)\oplus \theo (-l'-l)
\end{equation}
for some $l$ satisfying $(-l'/2) \leq l \leq m$ (here, as before, $f$ denotes the generic fiber of the projection $\pi:S\rightarrow C$).  We shall prove that if one of the Bia{\l}ynicki-Birula strata of $\psstack$ for the action of $\cs$ has nonempty intersection with $\psstack_m$ then that stratum is contained in $\psstack_m$.  We will next show that each $\psstack_m$ is quasiprojective (in particular, is K\"{a}hler) and consequently admits an ordering of its strata of the desired type; then, by semicontinuity of the ``generic fiber type'' $l$ in Equation \eqref{splittype}, no stratum in $\psstack_{m+1} \setminus \psstack_m$ lies above any stratum in $\psstack_m$ in the ordering of strata of $\psstack_{m+1}$, and the localization argument applies to all $\psstack_m$.

We begin, then, by showing that $l$ is constant on strata of the decomposition for the action of $\cs$ on $\psstack$.
\begin{prop}
\mbox{}
\begin{enumerate}
\item Suppose that $(\cE,\phi)\in\psstack(\spec\boldc)$, and that
\bd 
\cE\big|_f \cong \theo(l)\oplus\theo(-l'-l).
\ed
Then, writing $\ds (\cE',\phi') = \lim_{\lambda\rightarrow 0}\lambda\cdot(\cE,\phi)$, one also has
\bd
\cE'\big|_f \cong \theo(l)\oplus \theo(-l'-l).
\ed
\item If $(\cE,\phi)$ and $(\cE',\phi')$ lie in the same component of the fixed-point set of $\cs$ in $\psstack$, then $\cE\big|_f \cong \cE'\big|_f$.
\end{enumerate}
\end{prop}
\begin{proof}
(1) Suppose $U\subseteq C$ is an open set and $V=\pi\inv(U)\subseteq S$, for which
\bd
\cE\big|_V \cong \theo (l\sigma)\big|_V \oplus \theo((-l'-l)\sigma)\big|_V.
\ed
Then by the construction of limits for orbits in Corollary \ref{corofextlemma}, one has
\begin{align*}
\lim_{\lambda\rightarrow 0}\lambda\cdot(\cE\big|_V) &\cong \lim_{\lambda\rightarrow 0}\lambda\cdot\theo(l\sigma)\big|_V \oplus \lim_{\lambda\rightarrow 0}\lambda\cdot \theo((-l'-l)\sigma)\big|_V\\
& \cong \theo(l\sigma)\big|_V \oplus \theo((-l'-l)\sigma)\big|_V
\end{align*}
as desired. \newline
 
(2) Since the components of the fixed-point set of $\cs$ in $\psstack$ are smooth, the conclusion follows, provided any $\cs$-equivariant deformation of  $\theo_{\pline} (l)\oplus\linebreak\theo_{\pline}(-l'-l)$ parametrized by an Artinian local $\boldc$-algebra is trivial; this, in turn, follows if
\bd
\Ext^1_{\pline}\bigl(\theo(l)\oplus\theo(-l'-l),\theo(l)\oplus\theo(-l'-l)\bigr)^{\cs} = 0.
\ed
This $\Ext^1$ reduces to $H^1\bigl(\pline,\theo(l'+2l)\oplus\theo(-l'-2l)\bigr)^{\cs}$; however the \v{C}ech complex computing $H^*\bigl(\pline,\sum_k\theo(k)\bigr)$ has no $\cs$-invariants in $H^1$: this is a consequence of the description in \cite{MR57:3116}, Theorem III.5.1. \end{proof}

It remains just to show that each $\psstack_m$ is quasiprojective; but the method of proof used in Theorem \ref{polystableisgood} to show that $\psstack$ is a scheme will give the desired result, provided we can show that the sheaves $\cE$ occurring in pairs $(\cE,\phi)$ that lie in $\psstack_m$ (with fixed Chern classes $c_1$ and $c_2$!) form a bounded family.

Suppose first that $\cE$ is of type {\em U}.  Then $\cE$ occurs in the canonical exact sequence
\bd
0\rightarrow M_1\otimes I_1 \rightarrow \cE \rightarrow M_2\otimes I_2\rightarrow 0.
\ed
  Write $M_1=\pi^*\cL\otimes\theo(l\sigma)$ and $M_2=\pi^*\bigl(E/\cL\bigr)\otimes\theo(-l'-l)\sigma)$ for some $l$ and some line subbundle $\cL\subset E$ (recall that $E=L_1\oplus L_2)$.  Then, as in Equation \eqref{lastone}, 
\begin{equation}\label{newc2}
c_2\cE = -(l'+2l)\deg \cL + c_2I_1 + c_2I_2.
\end{equation}
One has $l\geq -l'/2$ and so 
\bd
-(l'+2l)\deg \cL \geq 0.
\ed
Moreover, since $\cE$ is of type {\em U}, one obtains $-(l'+2l)<0$, and consequently for each $l$ satisfying $(-l'/2) < l \leq m$ there are only finitely many possible values of $\deg\cL\leq 0$ for which is it possible to satisfy Equation \eqref{newc2} with $c_2I_1, c_2I_2 \geq 0$.  But now for each of the finitely many possible collections of fixed data $l$, $\deg\cL$, $c_2I_1$ and $c_2I_2$, such 1-extensions are parametrized by a subscheme of a projective bundle over 
\bd
\operatorname{Quot}_{E/C}(\deg\cL)\times \Hilb_S^{c_2I_1} \times \Hilb_S^{c_2I_2}
\ed
(here $\operatorname{Quot}_{E/C}(\deg\cL)$ denotes the Quot-scheme parametrizing subsheaves of $E$ of degree $\deg\cL$), thus giving boundedness.

Finally, suppose $\cE$ is of type {\em E}.  After twisting with $\theo(l'\sigma)$, we may assume $l'=0$.  If $\cE$ is not locally free, we have
\bd
c_2\cE = c_2\cE^{**} + \text{length}\bigl(\cE^{**}/\cE\bigr);
\ed
moreover the formula of Section 4 of \cite{MR85i:14009}  shows that if $\cE^{**}$ has canonical exact sequence
\begin{equation}\label{ddcanon}
0\rightarrow \cF' \rightarrow \cE^{**} \rightarrow \cI_{Z\subset Y} \rightarrow 0,
\end{equation}
then $c_2\cE^{**} = \deg Z$, so $c_2\cE^{**}\geq 0$ and there are only finitely many possible choices for $c_2\cE^{**}$ and $\text{length}\bigl(\cE^{**}/\cE\bigr)$.  Since a relative Quot-scheme for a bounded family (here we have in mind the Quot-scheme that parametrizes the quotients $\cE^{**}\rightarrow \cE^{**}/\cE$ of the vector bundle $\cE^{**}$) is itself a scheme of finite type, it will be enough to show that for each fixed choice of $c_2\cE^{**}$, the family of such vector bundles $\cE^{**}$ with the prescribed second Chern class is bounded.

Recall that $\cF' = \pi^*E'$ for some rank two subsheaf $E'\subseteq E$.  Restricting the exact sequence \eqref{ddcanon} to $D$ gives
\bd
\cF'\big|_D \rightarrow \cE^{**}\big|_D \rightarrow \cI_{Z\subset Y}\big|_D \rightarrow 0,
\ed
which must in fact be left exact as well (because $\cF'\big|_D \rightarrow \cE^{**}\big|_D$ is generically injective and $\cF'\big|_D$ is torsion-free).  There is a short exact sequence
\bd
0\rightarrow \cI_{Z\subset Y} \rightarrow \theo_Y \rightarrow \theo_Z \rightarrow 0
\ed
coming from the isomorphism $\cI_{Z\subset Y} \cong \cI_Z/\cI_Y$, 
so the length of $I_{Z\subset Y}\big|_D$ satisfies
\begin{align*}
\text{length}\Bigl(\cI_{Z\subset Y}\big|_D\Bigr) & \leq \text{length}\Bigl(\theo_Y\big|_D\Bigr) + \text{length}\Bigl(\operatorname{Tor}^1\bigl(\theo_D,\theo_Z\bigr)\Bigr)\\
& \leq 2\,\text{length}\bigl(\theo_Z\bigr).
\end{align*}
In particular, 
\bd
\deg E -2\,\text{length}\bigl( \theo_Z\bigr) \leq \deg E'\leq \deg E.
\ed
As a result, the collection of bundles $\cE^{**}$ may be parametrized by subschemes of projective bundles over finitely many products $\operatorname{Quot}\times\Hilb$ of a Quot-scheme for subsheaves of $E$ and a Hilbert scheme of points on $S$.  This completes the proof of boundedness of the collection of sheaves $\cE$ appearing in pairs $(\cE,\phi)\in\M_S(E)_m(\spec\boldc) = \psstack_m(\spec\boldc).$

\subsection{Proof of Theorem \ref{homologybasisthm}}

The description of Equation $\ref{fpeq}$ makes it clear what the components of the fixed-point set of  $\psstack$ must be: they are isomorphic to products 
\bd
\Hilb_L^{\alpha} \times \Hilb _L^{\beta}
\ed
of $\cs$-invariant sets of Hilbert schemes of points in the total space of $L$; these may be described (see \cite{MR1711344}) by partitions $\alpha$, $\beta$.  We need therefore only compute the rank of the negative normal bundle along each component $\Hilb_L^{\alpha} \times \Hilb _L^{\beta}$.  We may compute this rank at the general point of the component, that is, one for which
\bd
\supp (\theo/I_1) \cap \supp (\theo/I_2) = \emptyset.
\ed

Suppose $(\cE,\phi)$ is such a general point; then the tangent space to $\psstack$ at $(\cE,\phi)$ is $\Ext^1\bigl(\cE,\cE(-D)\bigr)$.  Now $\Ext^1(\cE,\cE)$ splits into four components,
\begin{align}
 & \Ext^1\bigl(\pi^*L_1 \otimes \theo (l\sigma)\otimes I_1, \pi^*L_1 \otimes \theo (l\sigma)\otimes I_1\bigr),\\
 & \Ext^1\bigl(\pi^*L_2\otimes \theo ((-l'-l)\sigma)\otimes I_2,\pi^*L_2\otimes \theo ((-l'-l)\sigma)\otimes I_2\bigr),\\
 & \Ext^1\bigl(\pi^*L_1\otimes \theo (l\sigma)\otimes I_1,\pi^*L_2\otimes \theo ((-l'-l)\sigma)\otimes I_2\bigr),\\
\intertext{and}
 & \Ext^1\bigl(\pi^*L_2\otimes \theo ((-l'-l)\sigma)\otimes I_2, \pi^*L_1\otimes \theo (l\sigma)\otimes I_1\bigr).
\end{align}

The first two factors are exactly the tangent spaces to the moduli spaces of rank one sheaves at $\pi^*L_1\otimes\theo (l\sigma)\otimes I_1$ and $\pi^*L_2\otimes\theo((l'-l)\sigma)\otimes I_2$, respectively, and the first-order deformations of these sheaves that vanish along $D$ are exactly the tangent spaces to the Hilbert schemes $T_{I_1}\Hilb_L^{\alpha}$ and $T_{I_2}\Hilb_L^{\beta}$, respectively.  The negative normal bundles of these schemes have ranks $|\alpha | - \ell (\alpha)$ and $|\beta | - \ell (\beta)$, respectively (see Lemma 7.6 of \cite{MR1711344}).  

It remains, then, to compute the negative weight spaces of the two remaining direct factors.  
We will use the short exact sequences

\bd
0\rightarrow \pi^*L_1\otimes \theo (l\sigma)\otimes I_1 \rightarrow \pi^*L_1\otimes\theo(l\sigma)
\rightarrow \pi^*L_1\otimes\theo (l\sigma)\otimes \theo / I_1 \rightarrow 0
\ed

and
\begin{multline*}
0\rightarrow \pi^*L_2\otimes \theo ((-l'-l)\sigma)\otimes I_2 \rightarrow \pi^*L_2\otimes\theo ((-l'-l)\sigma)\rightarrow\\
 \pi^*L_2\otimes\theo ((-l'-l)\sigma)\otimes \theo/I_2 \rightarrow 0
\end{multline*}
to induce long exact sequences in the first variable.  

Now
\bd
\Ext^1 \bigl(\pi^*L_1\otimes\theo (l\sigma)\otimes \theo/I_1, \pi^*L_2 \otimes\theo ((-l'-l)\sigma)\otimes I_2\bigr) = 0
\ed
and
\bd
\Ext^1 \bigl(\pi^*L_2\otimes\theo ((-l'-l)\sigma)\otimes \theo/I_2,  \pi^*L_1\otimes \theo (l\sigma)\otimes I_1\bigr) = 0
\ed
because the cosupports of $I_1$ and $I_2$ are disjoint.  Similarly, the $\Ext^2$ term
\bd
\Ext^2 \bigl(\pi^*L_1\otimes\theo(l\sigma), \pi^*L_2\otimes \theo ((-l'-l)\sigma)\otimes I_2\bigr)
\ed
and its counterpart vanish because they sit in the exact sequence
\begin{multline*}
\Ext^1 \bigl(\pi^*L_1\otimes\theo(l\sigma), \pi^*L_2\otimes\theo ((-l'-l)\sigma)\otimes\theo/I_2
\bigr)\rightarrow\\
 \Ext^2 \bigl(\pi^*L_1\otimes\theo(l\sigma), \pi^*L_2\otimes \theo ((-l'-l)\sigma)\otimes I_2\bigr)
\rightarrow\\
 \Ext^2 \bigl(\pi^*L_1\otimes\theo (l\sigma), \pi^*L_2\otimes\theo ((-l'-l)\sigma)\bigr)
\end{multline*}
and its counterpart, respectively, in which the left and right terms vanish: the left vanishes because it computes $H^1$ of a skyscraper sheaf, the right because it computes 
\bd
H^1\left(C, (L_1\inv\otimes L_2)^{\pm 1}\otimes \mathbf{R}^1\pi_*\theo (\pm(l'+2l)\sigma)\right),
\ed
which vanishes as in the proof of Theorem \ref{polystableisgood}.

Thus, to compute the remaining negative weight spaces of the tangent space to $\psstack$, it is enough to compute the negative weight spaces of
\begin{align}
 & \Ext^1\bigl(\pi^*L_1\otimes\theo(l\sigma),\pi^*L_2\otimes\theo ((-l'-l)\sigma)\otimes I_2\bigr),\\\label{ext1s}
 & \Ext^1\bigl(\pi^*L_2\otimes\theo ((-l'-l)\sigma),  \pi^*L_1\otimes\theo(l\sigma)\otimes I_1\bigr),\\\label{ext1s2}
& \Ext^2\bigl(\pi^*L_1\otimes\theo(l\sigma)\otimes\theo/I_1, \pi^*L_2\otimes\theo ((-l'-l)\sigma)\bigr),\\
\intertext{and}
& \Ext^2\bigl(\pi^*L_2\otimes\theo ((-l'-l)\sigma)\otimes\theo/I_2,  \pi^*L_1\otimes\theo(l\sigma)\bigr).
\end{align}

Note that we may (as we have done) omit the ideals in the second variables of the $\Ext^2$ groups here because of the assumption that the cosupports of $I_1$ and $I_2$ are disjoint.

We have exact sequences for the $\Ext^1$ groups of Equations \ref{ext1s} and \ref{ext1s2} in the second variable, but this time the vanishing of cohomology groups involving terms
\bd
(L_1\inv\otimes L_2)^{\pm 1} \otimes L^j
\ed
implies that in fact
\begin{multline}
\Ext^1\bigl(\pi^*L_1\otimes\theo(l\sigma),\pi^*L_2\otimes\theo ((-l'-l)\sigma)\otimes I_2\bigr)\cong\\
 \Ext^0 \bigl(\pi^*L_1\otimes\theo(l\sigma), \pi^*L_2\otimes\theo((-l'-l)\sigma)\otimes\theo/I_2
\bigr)
\end{multline}
and
\begin{multline}
\Ext^1\bigl(\pi^*L_2\otimes\theo ((-l'-l)\sigma),  \pi^*L_1\otimes\theo(l\sigma)\otimes I_1\bigr)\cong\\
 \Ext^0\bigl( \pi^*L_2\otimes\theo ((-l'-l)\sigma),  \pi^*L_1\otimes\theo(l\sigma)\otimes\theo/I_1\bigr).
\end{multline}
We now choose a sequence of integers $(w_1,w_2,w_3)\in\mathbf{Z}^3$, satisfying
\bd
w_3 \gg w_2 \gg w_1 > 0,
\ed
 and adjust the weights of the \mbox{$(\cs)^2\times\cs$} action on $\psstack$ accordingly by precomposing with the map \mbox{$(\cs)^2\times\cs \rightarrow (\cs)^2\times\cs$} determined by these three weights.  

We have reduced our computation to the computation of the negative weight spaces in the four groups
\begin{align}
& H^0 \bigl(\theo ((-l'-2l)\sigma)\otimes\theo/I_2\bigr) \underset{\boldc}{\otimes}\boldc (w_2-w_1),\\
& H^0 \bigl(\theo ((l'+2l)\sigma)\otimes\theo/I_1\bigr) \underset{\boldc}{\otimes}\boldc (w_1-w_2),\\
& \Ext^2 \bigl( \theo/I_1,\theo((-l'-2l)\sigma)\bigr)\underset{\boldc}{\otimes}\boldc (w_2-w_1),\\
\intertext{and}
& \Ext^2 \bigl(\theo/I_2,\theo((l'+2l)\sigma)\bigr)\underset{\boldc}{\otimes}\boldc (w_1-w_2).
\end{align}
Here $\boldc (\chi)$ denotes the one-dimensional representation of $\cs$ with weight $\chi$.

Moreover, we have
\begin{multline}
 H^0 \bigl(\theo ((-l'-2l)\sigma)\otimes\theo /I_2\bigr) \otimes\boldc (w_2-w_1)\cong\\
 H^0 (\theo /I_2)\otimes \boldc \bigl(w_2-w_1-w_3(l'+2l)\bigr)
\end{multline}
and
\begin{multline}
H^0 \bigl(\theo ((l'+2l)\sigma)\otimes\theo/I_1\bigr) \otimes\boldc (w_1-w_2)\cong\\
 H^0 (\theo/I_1)\otimes \boldc \bigl(w_1-w_2+w_3(l'+2l)\bigr).
\end{multline}

Now, as in Lemma \ref{i2istrivial}, if we decompose $\theo/I_j$ by weight, say
\bd
\theo/I_j = \bigoplus_{n\geq 0} V(-nw_3),
\ed
with filtered pieces
\bd
F(-nw_3) = \bigoplus_{m\geq n} V(-mw_3),
\ed
then $\cs$ acts in 
\bd
\Ext^2 \bigl(V(-nw_3),\theo\bigr) := \Ext^2 \bigl(F(-nw_3)/F(-(n+1)w_3), \theo\bigr)
\ed
with weight $(n+1)w_3$:  for, using the resolution
\bd
\theo(-\sigma-f) \rightarrow \theo(-f)\oplus\theo(-\sigma) \rightarrow \theo \rightarrow \boldc
\ed
gives
\bd
\Ext^2 (\boldc, \theo) = \operatorname{coker} \Bigl(\theo (\sigma)\oplus\theo (f)\rightarrow \theo(\sigma + f)\Bigr)
\ed
locally, and this latter group has local generator $z\inv\left(l^*/s^*\right)\inv$, where $z$ is a local parameter for the curve $C$ and $l^*/s^*$ is a local variable in the fiber direction; but on this generator $\cs$ acts (under the inverse action!) with weight $w_3$.

We may summarize the weight space decompositions, then, as follows.  Write the partitions $\alpha$ and $\beta$ as
\bd
\alpha = 1^{a_1}2^{a_2}3^{a_3}\dots
\ed
and
\bd
\beta = 1^{b_1}2^{b_2}3^{b_3}\dots.
\ed
Then
\begin{multline}
H^0 (\theo/I_2)\otimes \boldc \bigl(w_2-w_1-w_3(l'+2l)\bigr)= \\
 \bigoplus_{j\geq 0} \bigoplus_{i=0}^{j-1} \boldc ^{b_j} \bigl(w_2-w_1-w_3(l'+2l+i)\bigr),
\end{multline}
\begin{multline}
H^0(\theo/I_1)\otimes \boldc \bigl(w_1-w_2+w_3(l'+2l)\bigr)= \\
 \bigoplus_{j\geq 0}\bigoplus_{i=0}^{j-1} \boldc ^{a_j}\bigl(w_1-w_2+w_3(l'+2l-i)\bigr),
\end{multline}
\begin{multline}
\Ext^2(\theo/I_1,\theo)\otimes\boldc \bigl(w_2-w_1-w_3(l'+2l)\bigr)= \\
 \bigoplus_{j\geq 0} \bigoplus_{i=0}^{j-1} \boldc ^{a_j}\bigl(w_2-w_1+w_3 (i+1 - (l'+2l))\bigr),
\end{multline}
and
\begin{multline}
\Ext^2(\theo/I_2,\theo)\otimes \boldc \bigl(w_1-w_2+w_3(l'+2l)\bigr)= \\
 \bigoplus_{j\geq 0}\bigoplus_{i=0}^{j-1}\boldc ^{b_j} \bigl(w_1-w_2+w_3(i+1 + (l'+2l))\bigr).
\end{multline}
Adding all terms, we get
\begin{align*}
\bigoplus_{j\geq 0}\bigoplus_{i=0}^{j-1} \biggl[ \boldc ^{b_j} & \bigl(w_2-w_1-w_3(l'+2l+i)\bigr)\oplus \boldc ^{b_j}\bigl(w_1-w_2+w_3(l'+2l+i+1)\bigr) \oplus \\
& \boldc ^{a_j}\bigl(w_2-w_1+w_3(i+1-l'-2l)\bigr)
\oplus \boldc ^{a_j}\bigl(w_1-w_2+w_3(l'+2l-i)\bigr)\biggr] .
\end{align*}
Now, under our assumption that 
\bd
w_3 \gg w_2 \gg w_1,
\ed
we find that
\begin{align}
w_2 -w_1 -w_3(l'+2l+i)<0 & \text{\hspace{2.5em}iff\hspace{.5em}} & l'+2l+i & \geq 1,\\
w_1-w_2+w_3(l'+2l+i+1)<0 & \text{\hspace{2.5em}iff\hspace{.5em}} & l'+ 2l + i & \leq -1,\\
w_2-w_1+w_3(-l'-2l+i+1)<0 & \text{\hspace{2.5em}iff\hspace{.5em}} & l'+2l-i & >1,\\
\intertext{and}
w_1-w_2+w_3(l'+2l-i) <0 & \text{\hspace{2.5em}iff\hspace{.5em}} & l'+2l-i & \leq 0.
\end{align}

In other words, for each choice of $i$ and $j$, we obtain in the negative weight space a copy of 
\begin{align*}
\boldc^{a_j}\oplus\boldc^{b_j} & \text{\hspace{3.5em}provided\hspace{.5em}} & l'+2l & \neq -i, i+1\\
\boldc^{a_j} & \text{\hspace{3.5em}provided\hspace{.5em}} & l'+2l & = -i\\
\boldc^{b_j} & \text{\hspace{3.5em}provided\hspace{.5em}} & l'+2l & = i+1.
\end{align*}

This may be summarized as follows: we obtain $\boldc^d$, where $d$ is given by
\bd
d = \sum_{j>i\geq 0} b_j (1 - \delta_{l'+2l+i,0}) + \sum_{j>i\geq 0} a_j (1 - \delta_{l'+2l-i-1,0}).
\ed
The total rank of the negative normal bundle (expressed in terms of the partitions $\alpha$ and $\beta$ and the degrees $l'$ and $l$) is then
\begin{align*}
d(\alpha, \beta, l, l') = \biggl[ |\alpha | - & \ell (\alpha) +\sum_{j>i\geq 0} a_j(1 - \delta_{l'+2l-i-1,0})  \biggr]\\
+ & \biggl[ |\beta | - \ell (\beta) + \sum_{j>i\geq 0}b_j (1 - \delta_{l'+2l+i,0}) \biggr] ,
\end{align*}
which gives the expression for Equation \ref{shiftindex}.

Equation \ref{decompositionformula} follows immediately from the fact that the Abel-Jacobi map for $C$ identifies $\Symp^{j+1} C$ as a $\boldp^j$-fibration over $C$.  This completes the proof of Theorem \ref{homologybasisthm}.\qedsymbol

\bibliographystyle{amsalpha}
\providecommand{\bysame}{\leavevmode\hbox to3em{\hrulefill}\thinspace}

\end{document}